\documentclass[11pt]{amsart}
\usepackage{amsmath, amssymb, amscd, amsthm, mathrsfs, url ,pinlabel, verbatim, lipsum, wrapfig}
\usepackage[text={6.5in,9.5in}, centering, letterpaper, dvips]{geometry}
\usepackage{color, multirow, dcpic, latexsym, pictexwd, graphicx, epstopdf, tikz-cd, wrapfig, float}
\usepackage{verbatim, relsize, latexsym, latexsym, tikz, mdframed, float, mathtools, flowchart, multicol, stmaryrd, caption, capt-of, old-arrows, dirtytalk, footmisc, setspace, enotez, tabularx, pifont, tensor, csquotes}

\usepackage[backref=page, linktocpage=true]{hyperref}

\renewcommand*{\backref}[1]{}
\renewcommand*{\backrefalt}[4]{%
    \ifcase #1 (Not cited.)%
    \or        (Cited on page~#2.)%
    \else      (Cited on pages~#2.)%
    \fi}

\definecolor{maroon}{rgb}{0.5, 0.0, 0.0}

\hypersetup{
  colorlinks   = true, 
  urlcolor     = maroon, 
  linkcolor    = maroon, 
  citecolor   = maroon 
}

\usetikzlibrary{positioning, arrows.meta}
\tikzstyle{bag} = [align=center]

\theoremstyle{definition}

\newtheorem{introthm}{Theorem}

\newtheorem{introprob}{Problem}

\newcommand\Z{\mathbb{Z}}
\newcommand\Q{\mathbb{Q}}
\newcommand\TZ{\Theta_{\mathbb{Z}}^3}
\newcommand\TQ{\Theta_{\mathbb{Q}}^3}
\newcommand\HC{\Theta^n}
\newcommand\HCZ{\Theta_{\mathbb{Z}}^n}
\newcommand\intB{\smash{\mathring{\nu(K)}}}
\newcommand\intC{\smash{\mathring{\nu(K_1)}}}
\newcommand\intD{\smash{\mathring{\nu(K_2)}}}
\newcommand\intE{\smash{\mathring{B^3}}}

\let\footnote=\endnote

\title{A Survey of the Homology Cobordism Group}

\author{O{\u{g}}uz \c{S}avk}
\address{Department of Mathematics, Stanford University, Stanford, CA 94305, USA and Department of Mathematics, Bo\u{g}az{\i}\c{c}{\i}  University, Bebek, Istanbul 34342, Turkey.}
\email{\url{oguzsavk@stanford.edu}, \url{oguz.savk@boun.edu.tr}}
\urladdr{\url{https://sites.google.com/view/oguzsavk/}}
\date{}

\begin{document}

\begin{abstract}

In this survey, we present most recent highlights from the study of the homology cobordism group, with a particular emphasis on its long-standing and rich history in the context of smooth manifolds. Further, we list various results on its algebraic structure and discuss its crucial role in the development of low-dimensional topology. Also, we share a series of open problems about the behavior of homology $3$-spheres and the structure of $\TZ$. Finally, we briefly discuss the knot concordance group $\mathcal{C}$ and the rational homology cobordism group $\TQ$, focusing on their algebraic structures, relating them to $\TZ$, and highlighting several open problems. The appendix is a compilation of several constructions and presentations of homology $3$-spheres introduced by Brieskorn, Dehn, Gordon, Seifert, Siebenmann, and Waldhausen.

\end{abstract}
\maketitle

\begin{spacing}{1.12}
\setcounter{tocdepth}{3}
\tableofcontents
\end{spacing}

\vspace{-1.5 em}

\section{A Promenade around Smooth Manifolds}

All $n$-dimensional manifolds ($n$-manifolds for short) with or without boundaries are chosen to be compact, connected, oriented, and smooth. Otherwise, the type of the manifold is specified. The boundary of a manifold $M$ is denoted by $\partial M$, and $-M$ stands for $M$ with the opposite orientation. The connected sum operation between two manifolds is denoted by $\#$. A diffeomorphism (resp. homeomorphism, and piecewise linear homeomorphism) indicates a smooth (resp. continuous, and continuous and piecewise linear) bijective map between manifolds with a smooth (resp. continuous, and continuous and piecewise linear) inverse.

\subsection{The Predecessor: \texorpdfstring{$\HC$}{Tn}}

An $n$-manifold $M$ with $\partial M  = \emptyset$ is called a \emph{homotopy $n$-sphere} if $M$ has the same homotopy type as the unit $n$-dimensional sphere $S^n$, i.e., $M \simeq S^n$.  The \emph{$n$-dimensional homotopy cobordism group} $\HC$ is defined as \begin{center}
$\HC = \{ \text{homotopy} \ n\text{-spheres up to diffeomorphism} \} / \sim$ \end{center} where the equivalence relation \emph{h-cobordism}\footnote{The terms \say{h-cobordism} and \say{J-equivalence} were used interchangeably in these references.} $\sim$ is given for two arbitrary homotopy $n$-spheres $M_0$ and $M_1$ as 

\begin{minipage}{0.7\linewidth}

\[  M_0 \sim M_1 
 \iff  \begin{cases}
\ \text{There exists an} \ (n+1)\text{-manifold} \ W \ \text{such that} \\ 
   \bullet \ \partial W = -(M_0) \cup M_1, \\
  \bullet \ \text{The inclusions induce homotopy equivalences:} \\
\ \ \ \ \ \    M_0 \hookrightarrow W \hookleftarrow M_1 \ \ \Rightarrow \ \ M_0 \simeq W \simeq M_1.
  \end{cases} \]
  
\end{minipage}  
\begin{minipage}{0.3\linewidth}  
  
\includegraphics[width=0.87\textwidth]{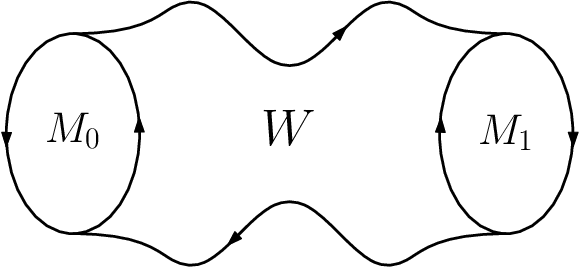}   
\end{minipage}  

\vspace{0.5 em}
  
After Milnor detected exotic $7$-spheres ($7$-manifolds homeomorphic but not diffeomorphic to $S^7$) in his groundbreaking work \cite{M56}, he also introduced the notion $\HC$ to study homotopy $n$-spheres in an unpublished note \cite{M59} and obtained some partial results on the orders of $\HC$. It forms an abelian group under the addition induced by connected sum. The zero element of $\HC$ is the homotopy cobordism class of $S^n$, and the inverse elements come with opposite orientation. Later, Kervaire and Milnor elaborated the structure of $\HC$ systematically in their celebrated article \say{\emph{Groups of homotopy spheres: I}} \cite{KM63}. 

Kervaire and Milnor were able to prove the following powerful statement, independent of the seminal articles of Connell \cite{C67}, Newman \cite{N66}, Smale \cite{S61}, Stallings \cite{St60}, and Zeeman \cite{Z61} about the topological Poincar\'e conjecture and the piecewise linear Poincar\'e conjecture in higher dimensions.\footnote{The topological (resp. piecewise linear, and smooth) Poincar\'e conjecture asserts that every topological (resp. piecewise linear, and smooth) homotopy $n$-sphere is homeomorphic (resp. piecewise linear homeomorphic, and diffeomorphic) to $S^n$. The topological and piecewise linear Poincar\'e conjectures were both proved for $n \geq 5$ in aforementioned articles. The particular case of $n=4$ for the topological Poincar\'e conjecture was shown in the seminal article of Freedman \cite{F82}, also see the book of Behrens, Kalm\'ar, Kim, Powell, and Ray \cite{BKKPR21}. The piecewise linear Poincar\'e conjecture in dimension 4 is still an open problem and is equivalent to the smooth Poincar\'e conjecture in dimension $4$ as a result of the articles of Cerf \cite{Ce68} and Hirsch and Mazur \cite{HM74}, see Rudyak's books \cite[IV.4.27(iv)]{Ru98} and \cite[6.7~Remark]{Ru16} for a detailed explanation.} Furthermore, they created the famous table with a single unknown value, depicted in Table~\ref{tab:elts}.

\begin{introthm}[Theorem~1.2, \cite{KM63}]
For $n\neq 3$, the group $\HC$ is finite.
\end{introthm}

\begin{table}[htbp]
\centering
\renewcommand{\arraystretch}{1.2}
\begin{tabular}{|l|l|l|l|l|l|l|l|l|l|l|l|l|l|l|l|l|l|l|}
\hline
$n$ & $1$ & $2$ & $3$ & $4$ & $5$ & $6$ & $7$ & $8$ & $9$ & $10$ & $11$ & $12$ & $13$ & $14$ & $15$ & $16$ & $17$ & $18$ \\ \hline
$\vert \HC \vert$  & $1$ & $1$ & $?$ & $1$ & $1$ & $1$ & $28$ & $2$ & $8$ & $6$ & $992$ & $1$ & $3$ & $2$ & $16256$ & $2$ & $16$ & $16$ \\ \hline
\end{tabular}
\vspace{0.5 em}
\caption{The orders of $\HC$ for $1\leq n \leq 18$.}
\label{tab:elts}
\vspace*{-2 em}
\end{table}

The classical results of Moise \cite{M52a, M52b} showed that every topological $3$-manifold has a unique smooth structure. After the confirmation of the last topological Poincar\'e conjecture, the missing point in Table~\ref{tab:elts} was clarified as an immediate consequence of Perelman's breakthrough.

\begin{introthm}[\cite{P1,P3,P2}]
The group $\Theta^3$ is trivial, hence $\vert \Theta^3 \vert = 1$.
\end{introthm}

Kervaire and Milnor never published \say{\emph{Groups of homotopy spheres: II}}; however, Levine's lecture notes \cite{L85} can be considered as its sequel paper.\footnote{See the introduction of \cite{L85}. Also consult Milnor's survey \cite[Pg.~805]{M11}, and the commentary of Ranicki and Webber on the correspondence of Kervaire and Milnor around the 1960s \cite{RW15}.} Finding the order of $\HC$ for each value of $n$ is a very challenging problem in algebraic and geometric topology. Moreover, it is closely tied to the smooth Poincar\'e conjecture in higher dimensions.\footnote{The smooth Poincar\'e conjecture is false in general. For precise expositions, consult the introduction of \cite{WX17} and also see the papers of Isaksen \cite{I19} and Isaksen, Wang, and Xu \cite{IWX20}.} For the state of the art regarding the order of $\HC$, one can see \cite[Table~1]{IWX20b}. 

Further discussions and results about homotopy theoretical approaches to study $\Theta^n$ can be seen in excellent papers of Hill, Hopkins, and Ranevel \cite{ HHR16}, Wang and Xu \cite{WX17}, and Behrens, Hill, Hopkins, and Mahowald \cite{BHHM20}. 

\subsection{The Successor: \texorpdfstring{$\HCZ$}{TnZ}}

In a similar vein, a \emph{homology $n$-sphere} is an $n$-manifold $M$ with $\partial M = \emptyset$ such that $M$ has the same homology groups of $S^n$ in integer coefficients, i.e., $H_*(M;\Z) = H_*(S^n;\Z)$. The \emph{$n$-dimensional homology cobordism group} $\HC$ is formed as \begin{center}
$\HCZ = \{ \text{homology} \ n\text{-spheres up to diffeomorphism} \} / \sim_\Z$ \end{center} where the equivalence relation \emph{homology cobordism} $\sim_\Z$ is depicted for two arbitrary homology $n$-spheres $M_0$ and $M_1$ by \[ M_0 \sim_\Z M_1 
 \iff  \begin{cases}
\ \text{There exists an} \ (n+1)\text{-manifold} \ W \ \text{such that} \\  
 \bullet \ \partial W = -(M_0) \cup M_1, \\
  \bullet \ \text{The inclusions induce isomorphisms on all homology groups:} \\ \ \ \ M_0 \hookrightarrow W \hookleftarrow M_1 \ \ \Rightarrow \ \ H_*(M_0; \Z) \cong H_*(W; \Z) \cong H_*(M_1; \Z). \end{cases} \]  
  
Inspired by the novel work of Kervaire and Milnor, Gonz\'{a}lez-Acu\~{n}a defined the object $\HCZ$ to decipher the homology $n$-spheres in his Ph.D. thesis \say{\emph{On homology spheres}} \cite{GA70}. Similarly, $\HCZ$ admits an abelian group structure with the summation induced by connected sum. The homology cobordism class of $S^n$ serves as the identity element of $\HCZ$. Besides, inverse elements can be obtained by reversing the orientation. 

Using surgery theory and Milnor's $\pi$-manifolds,\footnote{Similarly, \say{$\pi$-manifold} and \say{$s$-parallelizable} as well as \say{surgery} and \say{spherical modification} were different names for the same notion. An $n$-manifold $M \subset \mathbb{R}^{n+q}$ is called a \emph{$\pi$-manifold} if its normal bundle $\nu(M)$ is trivial, i.e., $\nu(M)$ is diffeomorphic to $M \times \mathbb{R}^q$.} Gonz\'{a}lez-Acu\~{n}a was able to construct a group isomorphism between $\HC$ and $\HCZ$ unless $n=3$. Hence, they are algebraically identical except for the single case of $n=3$.

\begin{introthm}[Theorem~I.2, \cite{GA70}]
\label{Gonz}
For $n\neq 3$, $\HCZ$ is isomorphic to $\HC$. Therefore, $\HCZ$ is finite unless $n=3$.
\end{introthm}

It should be very interesting to compare Gonz\'{a}lez-Acu\~{n}a's elegant theorem with the following achievement of Kervaire which was published around the same time.

\begin{introthm}[Theorem~3, \cite{K69}]
For $n\geq 5$, let $M$ be a homology $n$-sphere. Then there exists
a unique homotopy sphere $\Sigma_M$ such that $M \# \Sigma_M$ bounds a contractible $(n+1)$-manifold.
\end{introthm}

\subsection{The Aberrant: \texorpdfstring{$\TZ$}{T3Z}}

The isomorphism of Gonz\'{a}lez-Acu\~{n}a cannot be valid for the last case $n=3$ due to the famous invariant of Rokhlin \cite{R52}. There is a surjective group homomorphism from the $3$-dimensional homology cobordism group (the homology cobordism group for short) to the cyclic group of order two \begin{center} $\mu: \TZ \to \Z_2, \ \ \ \mu(Y) = \sigma(W) / 8 \mod 2$
\end{center} where $W$ is any $4$-manifold with a $\Z_2$-valued even intersection form,\footnote{For the other reformulations of the Rokhlin invariant $\mu$ in terms of the characterization of a $4$-manifold, see the recent ICM 2022 paper of Finashin, Kharlamov, and Viro \cite{FKV20}.} $\partial W = Y$, and $\sigma(W)$ denotes the signature of $W$. 

The homology cobordism invariance of the Rokhlin invariant $\mu$ was first observed in \cite[Section~I.5]{GA70}. See also \cite[Section~2]{G70} and \cite[Section~3.8]{FK20}. Since the Poincar\'e homology sphere $\Sigma(2,3,5)$ (see Section~\ref{appendix} for its several descriptions) uniquely bounds the negative-definite plumbing $-E_8$ of signature $-8$, we have $\mu(\Sigma(2,3,5)) = 1$. Therefore, it is not homology cobordant to $S^3$, and we conclude:

\begin{introthm}[\cite{R52}; Section~I.5, \cite{GA70}]
The group $\TZ$ is non-trivial.\footnote{Note that the homology cobordism group also appeared with notations $\Theta_3^H$ or $\mathscr{H}^3$ in the literature of the 1970-80s.}
\end{introthm}

The non-triviality of $\TZ$ is sensitive to both homology and smoothness conditions on the cobordism $4$-manifold. The group would be trivial if at least one of these conditions were removed. See the articles by Rokhlin \cite{R51} and Freedman \cite{F82} respectively. Also, $\TZ$ is countable by the classical results of Moise \cite{M52a, M52b}.

Until the 1980s, the only known invariant of $\TZ$ was the Rokhlin invariant $\mu$ and there was a belief that it might be an isomorphism. However, it later turned out that $\TZ$ is far from being finite. The understanding of the infinitude of $\TZ$ has led to the construction of numerous invariants of homology $3$-spheres.

The seminal work of Matumoto \cite{M78} and Galewski and Stern \cite{GS80} yielded a rich connection between the Rokhlin invariant $\mu$, the group $\TZ$, and the triangulation conjecture. Manolescu revolutionized low-dimensional topology by introducing the Seiberg-Witten (monopole) $\mathrm{Pin}(2)$-equivariant Floer homology, constructing the $\beta$-invariant, and disproving the triangulation conjecture \cite{M16}. His $\beta$-invariant is an integer lift of the Rohklin invariant $\mu$ and its existence rejects the triangulation conjecture by relying on the articles \cite{M78, GS80}. Consult Section~\ref{torsion} for more details. The several variations of Manolescu's Floer homotopic approach have led to the invention of new powerful theories and sensitive invariants of knots and manifolds. Recently, there has also been increased activity in studying $\TZ$ using techniques from $SU(2)$-gauge theory, following the work of Daemi \cite{Dae20}.

\vspace{0.25 em}

\noindent \begin{minipage}{0.83\textwidth}

\hspace{1.2 em} Homology cobordism is closely related to the concepts of \emph{knot concordance} and \emph{rational homology cobordism}, and both give rise to abelian groups $\mathcal{C}$ and $\TQ$, similar to $\TZ$. By the classical work of Gonz\'{a}lez-Acu\~{n}a \cite{G70}, Gordon \cite{G75}, and Casson and Gordon \cite{CG78}, there are natural mappings between these three abelian groups given by $(1/n)$-surgery on knots in the $3$-sphere $S^3_{1/n} (K)$ for any integer $n$, $p^r$-fold cyclic branched coverings of the $3$-sphere along knots $\Sigma_{p^r} (K)$ for any prime $p$ and $r \geq 1$, and inclusion $\psi$. Consult Section~\ref{Qcobordism}, Section~\ref{Cconcordance}, and Section~\ref{appendix} for further details.

\end{minipage}
\begin{minipage}{0.2\textwidth}
\vspace*{- 0.9 em}
\begin{center}
\begin{tikzcd}[row sep=tiny]
								& \TZ \arrow{dd}{\psi}			\\
  \mathcal{C} \arrow{ur}{S^3_{1/n}} \arrow{dr}{\Sigma_{p^r}} &              \\
										& \TQ
\end{tikzcd}
\end{center}

\end{minipage}

\vspace{0.5 em}

\noindent \begin{minipage}{0.65\linewidth}
\hspace{1.2 em} In a nutshell, we create this table to reflect the sharp contrast between the homology cobordism group $\TZ$ and all other homotopy and homology cobordism groups. One can access the most recent information about the orders of $\HC$ from the article of Isaksen, Wang, and Xu \cite{IWX20b}. 
\end{minipage}
\hspace{2 em}
\begin{minipage}{0.35\linewidth}
\vspace*{- 0.4 em}
\begin{tabular}{l|ll|}
\cline{2-3}
& \multicolumn{2}{l|}{\textbf{\ \ \ Order}}       \\ \hline
\multicolumn{1}{|l|}{\textbf{Dimension}} & \multicolumn{1}{l|}{\ $\HC$} & \multicolumn{1}{l|}{\ $\HCZ$}  \\ \hline
\multicolumn{1}{|l|}{\ \ \ $n \neq 3$} & \multicolumn{1}{l|}{$< \infty$} & \multicolumn{1}{l|}{$< \infty$}  \\ \hline
\multicolumn{1}{|l|}{\ \ \ $n=3$} & \multicolumn{1}{l|}{$=1$} & \multicolumn{1}{l|}{$= \infty$}  \\ \hline
\end{tabular}

\end{minipage}

\vspace{0.7 em}

From now on, we will aim to approach all results that arise around the homology cobordism group $\TZ$ from a broad, comprehensive, and historical perspective. Our additional purpose is to present various open problems of homology $3$-spheres in the context of homology cobordism. Finally, we will discuss the knot concordance group $\mathcal{C}$ and the rational homology cobordism group $\TQ$ by eleborating their most recent algebraic structure, relating them to $\TZ$, and posing several open problems. Most of the problems raised in this survey are well-known in the field in general. We hope that our efforts will have a positive impact, and motivate on readers to investigate and study the homology cobordism group $\TZ$ in the future.

\section{The Structure of \texorpdfstring{$\TZ$}{T3Z}}
\subsection{Subgroups and Summands of \texorpdfstring{$\TZ$}{T3Z}} 
The celebrated work of Donaldson was a cornerstone in the history of low-dimensional topology \cite{D83}. Motivated by his article, Fintushel and Stern studied the gauge theory of orbifolds, produced the gauge theoretical $R$-invariant for Seifert fibered homology spheres, and provided the first existence of an infinite subgroup in the homology cobordism group. 

\begin{introthm}[Theorem~1.2, \cite{FS85}]
\label{FS}
The group $\TZ$ has a $\Z$ subgroup generated by the Poincar\'e homology sphere $\Sigma(2,3,5)$.
\end{introthm}

The extended version of Donaldson's diagonalization theorem \cite{D87} recovers Theorem~\ref{FS} as follows: One can use $\Sigma(2,3,5)$ to construct a closed $4$-manifold whose non-diagonalizable intersection form is $n E_8$ for arbitrary value of $n$. This obstructs the existence of any homology cobordism between $S^3$ and a finite number of self connected sums of $\Sigma(2,3,5)$.

Converting the ideas on end-periodic $4$-manifolds in the work of Taubes \cite{T87} to cylindrical end $4$-manifolds and using the Fintushel-Stern $R$-invariant, Furuta showed the first existence of an infinitely generated subgroup \cite{F90}.

\begin{introthm}[Theorem~2.1, \cite{F90}]
\label{F}
The group $\TZ$ has a $\Z^\infty$ subgroup\footnote{In our convention, $\Z^\infty$ always stands for $\bigoplus_{n=1}^\infty \Z$.} in $\TZ$ generated by the family of Brieskorn spheres $\{ \Sigma(2,3,6n-1)\}_{n=1}^\infty$.
\end{introthm}

The eminent article of Floer \cite{F88} changed the flow of the history of low-dimensional topology dramatically. Given a homology $3$-sphere $Y$, his theory of instanton homology can be defined over the Yang-Mills equations on $Y \times \mathbb{R}$. This novel invariant is an infinite dimensional analogue of the Morse homology.

The next achievement about the algebraic structure of $\TZ$ was owed to Fr{\o}yshov \cite{Fro02}. His approach relied on the equivariant structure on Floer's instanton (Yang-Mills) homology and he constructed the $h$-invariant, a surjective group homomorphism $h: \TZ \to \Z$.

\begin{introthm}[Theorem~3, \cite{Fro02}]
\label{Fr}
The group $\TZ$ has a $\Z$ summand generated by the Poincar\'e homology sphere $\Sigma(2,3,5)$.
\end{introthm}

Ozsv\'ath and Szab\'o developed the theory of Heegaard Floer homology in a series of prominent articles \cite{OS03b, OS04, OS03c}. Since then it has been used to answer various problems in low-dimensional topology and several new versions emerged successively, see the comprehensive surveys of Ozsv\'ath and Szab\'o \cite{OS04e} and Juh\'asz \cite{J15}. Later, Hendricks and Manolescu introduced involutive Heegaard Floer homology \cite{HM17} and this new theory exploits the conjugation symmetry on a Heegaard Floer complex of the Heegaard Floer homology. Also, it is conjecturally a $\Z_4$-equivariant version of Seiberg-Witten $\mathrm{Pin}(2)$-equivariant Floer homology established by Manolescu \cite{M16}.

The most recent impressive progress about deciphering the algebraic complexity of the group $\TZ$ was achieved by Dai, Hom, Stoffregen, and Truong \cite{DHST18}. Using the machinery of involutive Heegaard Floer homology, they defined a new family of powerful and sensitive set of invariants $ \vec{f} = \{ f_k \}_{k \in \mathbb{N}}$, a surjective group homomorphism $\vec{f}: \TZ \to \Z^\infty$.\footnote{In \cite{R20}, Rostovtsev reinterpreted the homomorphisms of Dai, Hom, Stoffregen, and Truong by using the immersed curve machinery of Kotelskiy, Watson, and Zibrowius \cite{KWZ19}. In particular, he found a new epimorphism of $\TZ$ independent of  $\{ f_k \}_{k \in \mathbb{N}}$.}

\begin{introthm}[Theorem~1.1, \cite{DHST18}]
\label{DHST}
The group $\Theta^3_\mathbb{Z}$ has a $\Z^\infty$ summand generated by the family of Brieskorn spheres $\{\Sigma(2n+1,4n+1,4n+3) \}^\infty_{ n=1}$.
\end{introthm} 

Their proof subsumes several approaches and techniques that consecutively appeared in the literature of involutive Heegaard Floer homology \cite{HMZ18}, \cite{DM17}, \cite{DS17}, and \cite{HHL18}. Moreover, involutive Floer theoretic invariants have provided a major change for the understanding of the structure of $\TZ$ and its subgroups. For details of constructions and ideas, one can consult the survey of Hom \cite{H21}.

Relying on all these previous results, one may expect that there is no torsion part in the decomposition of $\TZ$, see Section~\ref{torsion} for details. In particular, Problem~\ref{isomorphism} and Problem~\ref{torprob} are complementary, and Problem~\ref{pqn} is a special case of Problem~\ref{isomorphism}. The author believes that the following problem will have a negative answer.

\begin{introprob}
\label{isomorphism}
Is $\TZ$ is isomorphic to $\Z^\infty$? 
\end{introprob}

Most instanton, Seiberg-Witten, and Floer theoretical invariants of homology $3$-spheres are sensitive to a preorder given by the negative-definite cobordisms. Thus, the further understanding of the structure of the homology cobordism group will be possible by realizing $\TZ$ as a partially-ordered group, rather than just a group. See, for instance, the recent work of Nozaki, Sato, and Taniguchi \cite[Section~1.3]{NST22}.

\begin{introprob}
Study the structure of $\TZ$ as an ordered group by forming filtrations, and completely describe subgroups and quotients.
\end{introprob}

\subsubsection{A Recovery: More about Subgroups of \texorpdfstring{$\TZ$}{T3Z}}
\label{recovery}

A $4$-manifold with boundary is called a \emph{homology $4$-ball} if it shares the same homology groups of the $4$-ball in integer coefficients. An easy algebraic topology argument indicates that a homology $3$-sphere is homology cobordant to $S^3$ if and only if it bounds a homology $4$-ball.

The Fintushel-Stern $R$-invariant leads to a powerful obstruction for homology $3$-spheres to bound homology $4$-balls, and hence contractible $4$-manifolds. It is easily computable due to the short-cut of Neumann and Zagier \cite{NZ85}. The non-zero values of the $R$-invariant provide the proofs of items \eqref{a} and \eqref{c} in Theorem~\ref{subgroup}. Further, these claims can be deduced by using the Ozsv\'ath-Szabo $d$-invariant \cite{OS03b}. See the papers of Tweedy \cite{T13} and Karakurt and the author \cite{KS20} for sample computations, which both depended on Floer homology of plumbings \cite{OS03a}, N\'emethi's lattice homology \cite{Nem05}, and lattice point counting technique of Can and Karakurt \cite{CK14}.\footnote{These three articles all provide equivalent but different descriptions of Heegaard Floer homology groups of Seifert fibered homology spheres.}

However, the item \eqref{b} in Theorem~\ref{subgroup} is a consequence of the non-vanishing of the Neumann-Siebenmann invariant $\bar{\mu}$ \cite{N80, S80}. The homology cobordism invariance of $\bar{\mu}$ for Seifert fibered homology spheres was first proved by Saveliev \cite{Sav98}, see also the paper of Dai and Stoffregen \cite{DS17} for a generalization of this result. Saveliev provided another proof for the item \eqref{b} in \cite{Sav98b} by using Furuta's $10/8 +2$ theorem \cite{Fu01}. Note that Furuta's result was a partial solution for Matsumoto's $11/8$ conjecture \cite{Mat82}. In this article, he also introduced a homology cobordism invariant called the \emph{bounding genus}. All other homology cobordism invariants behaved differently than $\bar{\mu}$ seem to vanish or not be arbitrarily large for this family, so they do not give further information about their homology cobordism classes. 

By following the work of Nozaki, Sato and Taniguchi \cite{NST22} and Baldwin and Sivek \cite{BS22}, the proofs of items \eqref{d} and \eqref{e} in Theorem~\ref{subgroup} can be deduced respectively. Moreover, the items \eqref{e2} and \eqref{e3} in Theorem~\ref{subgroup} are owed to the recent article of Daemi, Imori, Sato, Scaduto, and Taniguchi \cite{DISST22}. Note that the arguments of the latter two articles essentially require the result of the first one. Here, $\tau^\sharp$- and $\tilde{s}$-invariants are new instanton Floer theoeric invariants of knots \cite{BS22, DISST22}, and $h$ denotes the classical Fr{\o}yshov invariant which appeared in Theorem~\ref{Fr}, $\Gamma$ stands for the new invariant of knots, and both invariants are again derived from instanton Floer homology.

\begin{introthm}
\label{subgroup}
The following homology $3$-spheres individually generate $\Z$ subgroups in $\TZ$:
\begin{enumerate}
\item \label{a} $\Sigma(p,q,pqn-1)$ for each $n \geq 1$,
\item \label{b} $\Sigma(p,q,pqn+1)$ for each odd $n \geq 1$,\footnote{This result cannot be generalized to even values of $n$ since $\Sigma(2,3,13)$ and $\Sigma(2,3,25)$ are known to bound contractible $4$-manifolds.}
\item \label{c} $\Sigma(p_n, q_n, r_n)$ for each $n \geq 1$ where $p_n q_n + p_n r_n - q_n r_n = 1$,
\item \label{d} For each $n \geq 1$, $S^3_{1/n}(K)$ where $K$ is any knot in $S^3$ with $h(S^3_1(K)) < 0$, \footnote{Explicitly, the knot $K$ can be taken as the mirrors $K_n^*$ of the $2$-bridge knots $K_n$ corresponding to the rational number $\frac{2}{4n-1}$ as hyperbolic examples. For the satellite type of examples, one can pick the $(2,q)$-cable of any knot $K$ with odd $q \geq 3$, see \cite{NST22}.}
\item \label{e} For each $n \geq 1$, $S^3_{1/n}(K) $ where $K$ is any knot in $S^3$ with $\tau^\sharp (K) > 0$,\footnote{The knot $K$ can be chosen as either a knot having a transverse representative with positive self-linking number, or quasi-positive knot which is not smoothly slice, or an alternating knot with negative signature $\sigma$, under the convention $\sigma(T(2,3)) = -2$, see \cite{BS21} and \cite{BS22}.}
\item \label{e2} For each $n \geq 1$, $S^3_{1/n}(K) $ where $K$ is any knot in $S^3$ with $\tilde{s} (K) > 0$,\footnote{The knot $K$ can be chosen as either a quasi-positive knot which is not smoothly slice or an alternating knot with negative signature.}
\item \label{e3} For each $n \geq 1$, $S^3_{1/n}(K) $ where $K$ is any knot in $S^3$ with $\sigma(K) \leq 0$ and $\frac{1}{8} < \Gamma_K \left( -\frac{1}{2} \sigma(K) \right) $.\footnote{Under these conditions, Daemi, Imori, Sato, Scaduto, and Taniguchi provided two-parameter family of bridge knots $K_{m, n} = K(212mn-68n+53, 106m-34)$ ($m$ and $n$ are fixed) such that $(1/k)$-surgery on the mirrors of $K_{m, n}$ are linearly independent in the homology cobordism group yet $K_{m, n}$ are torsion in the algebraic concordance group of knots.}
\end{enumerate}
\end{introthm} 

Manolescu's invariants $\alpha, \beta, \gamma$ \cite{M16} and Hendricks-Manolescu's involutive $d$-invariants $\underline{d}, \overline{d}$ \cite{HM17}\footnote{Note that the involutive correction terms $\underline{d}$ and $\overline{d}$ in \cite{HM17} and Manolescu invariants $\alpha, \beta$ and $\gamma$ in \cite{M16} are not homomorphisms.} can be read off from the values of the Ozsv\'ath-Szab\'o $d$-invariant and the Neumann-Siebenmann $\bar{\mu}$-invariant, see articles Dai and Manolescu \cite{DM17} and Stoffregen, \cite{St20} for more details. In particular, the $R$-invariant of Fintushel-Stern \cite{FS85} is directly determined from a plumbing graph due to the short-cut of Neumann-Zagier \cite{NZ85}. Moreover, the $\bar{\mu}$-invariant of Seifert fibered homology spheres is same as the $w$-invariant of Fukumoto and Furuta \cite{FF00}, see the work of Fukumoto, Furuta, and Ue \cite{FFU01} and Saveliev \cite{Sav02c} for details. Therefore, we have the following several identities between homology cobordism invariants for a single Seifert fibered space $\Sigma = \Sigma(a_1, \ldots, a_n)$:

\begin{itemize}

\item $R\left( \Sigma\right) = -2e-3$,
\item $d\left( \Sigma\right) = \overline{d}\left( \Sigma \right)$,
\item $\bar{\mu}\left( \Sigma \right) = w \left( \Sigma \right) = -\frac{1}{2} \underline{d}\left( \Sigma \right) = -\beta\left( \Sigma \right) = -\gamma\left( \Sigma \right)$,
\item $\alpha\left( \Sigma \right) = \begin{cases} 
      \frac{1}{2} d\left( \Sigma \right), & \text{if} \ \frac{1}{2} d\left( \Sigma \right) = -\bar{\mu}\left( \Sigma \right) \mod 2, \\
      \frac{1}{2} d\left( \Sigma \right) + 1, & \text{otherwise}, 
\end{cases}$,
\item $\mu \left( \Sigma \right) = \bar{\mu}\left( \Sigma \right) = \alpha \left( \Sigma \right) = \beta \left( \Sigma \right) = \gamma \left( \Sigma \right)  \mod 2$.
\end{itemize}

After Furuta's work, the first recovery of the existence of $\Z^\infty$ subgroups of $\TZ$ was provided by Fintushel and Stern \cite[Theorem~5.1]{FS90} for the item \eqref{f} in Theorem~\ref{subgroups}. Their approach can be applied to item \eqref{g} in Theorem~\ref{subgroups} as well. These two results can be reproved successfully by using new gauge and instanton theoretic invariants of Daemi \cite{Dae20}, Nozaki, Sato and Taniguchi \cite{NST22}, and Baldwin and Sivek \cite{BS21, BS22}. However, the classical and involutive Heegaard Floer theoretical invariants cannot identify the linear independence of the item \eqref{f} in $\TZ$.

The Seiberg-Witten and/or Heegaaard Floer originated invariants may detect the linear independence of subfamilies of the item \eqref{g} in Theorem~\ref{subgroups}. In this regard, see the work of Stoffregen \cite{St17} and Dai and Manolescu \cite{DM17}. However, it is not easily doable in general, see the discussion in \cite{KS20} and \cite{KS22} and compare with \cite{St17} and \cite{DM17}.

For the proofs of items \eqref{h}, \eqref{i}, \eqref{i2}, and \eqref{i3} in Theorem~\ref{subgroups}, one can see the articles of Nozaki, Sato and Taniguchi \cite{NST22}, Baldwin and Sivek \cite{BS22}, and Daemi, Imori, Sato, Scaduto, and Taniguchi \cite{DISST22}. The methodology of \cite{NST22} and \cite{DISST22} both refer to the equivariant instanton Floer theory with Chern-Simons filtration, while \cite{BS21, BS22} uses the framed instanton homology. Notice that these articles all provide new invariants for homology $3$-spheres and knots.

\begin{introthm}
\label{subgroups}
The following infinite families of homology $3$-spheres generate $\Z^\infty$ subgroups in $\TZ$: 
\begin{enumerate}
\item \label{f} $\{ \Sigma(p,q,pqn-1)\}_{n=1}^\infty$,
\item \label{g} $\{ \Sigma(p_n, q_n, r_n)\}_{n=1}^\infty$ where $p_n q_n + p_n r_n - q_n r_n = 1$,
\item \label{h} $\{ S^3_{1/n}(K) \}_{n=1}^\infty$ for any knot $K$ in $S^3$ with $h(S^3_1(K)) < 0$,
\item \label{i} $\{ S^3_{1/n}(K) \}_{n=1}^\infty$ for any knot $K$ in $S^3$ with $\tau^\sharp (K) > 0$,\footnote{Since positive knots in $S^3$ are quasi-positive and not smoothly slice due to Rasmussen \cite{Ras10}, the work of Baldwin and Sivek also generalizes a result of Gompf and Cochran \cite{CG88}: $S^3_{1/n}(K)$ individually generates a $\Z$ subgroup in $\TZ$ when $K$ is a positive knot in $S^3$.}
\item \label{i2} $\{ S^3_{1/n}(K) \}_{n=1}^\infty$ for any knot $K$ in $S^3$ with $\tilde{s} (K) > 0$,
\item \label{i3} $\{ S^3_{1/n}(K) \}_{n=1}^\infty$ for any knot $K$ in $S^3$ with $\sigma(K) \leq 0$ and $\frac{1}{8} < \Gamma_K \left( -\frac{1}{2} \sigma(K) \right) $.
\end{enumerate}
\end{introthm} 

Since all current homology cobordism invariants are blind to detect the linear independence of $\{\Sigma(p,q,pqn+1)\}_{n=1, \text{odd}}^\infty$ in $\TZ$,  we curiously ask the following problem. On the other hand, these manifolds might be homology cobordant to each other in $\TZ$. If so, this will also be a very interesting result.

\begin{introprob}
\label{pqn}
Does the family $\{\Sigma(p,q,pqn+1)\}_{n=1, \text{odd}}^\infty$ generate a $\Z^\infty$ subgroup or a $\Z^\infty$ summand in $\TZ$?
\end{introprob} 

The $R$- and $w$-invariants were successfully generalized in the articles of Fintushel and Lawson \cite{FL86} and Fukumoto \cite{Fuk11} respectively. Given a Seifert fibered sphere $Y = \Sigma(a_1, \ldots ,a_n)$, we respectively denote these invariants by $R(Y,e)$ and $w(Y,m)$ and call the \emph{generalized $R$-invariant} and the \emph{generalized $w$-invariant} where $e$ is an integer depending on Euler number and some other constraints, and $m$ is a tuple of integers. The generalized $R$- and $w$-invariants are strictly more powerful than the classical $R$- and $w$-invariants, and provide more sensitive obstructions for the existence of homology cobordisms between homology $3$-spheres. In particular, a combinotorial formula for the generalized $R$-invariant was found by Lawson \cite{L87} so that $R\left( \Sigma, 1 \right) = R\left( \Sigma\right)$. For sample computations, see Fukumoto's article \cite[Section~6]{Fuk11}. Fukumoto also gave estimates for Matsumoto's bounding genera for homology $3$-spheres using $w$-invariants \cite{Fuk09}.

Using $\mathrm{Pin}(2)$-equivariant Seiberg-Witten Floer K-theory, Manolescu constructed the integer-valued homology cobordism invariant $\kappa$ \cite{Man14}. Recently, Ue proved that the behaviours of the $\kappa$ invariant and the minus version of the $\bar{\mu}$ invariant for Seifert fibered spheres are very similar \cite{U22}: $\kappa(Y) + \bar{\mu}(Y) = 0$ or $2$. Relying on the Seiberg-Witten Floer spectrum and $\mathrm{Pin}(2)$-equivariant KO-theory and inspiring the construction of the Manolescu $\kappa$-invariant, J. Lin extracted new invariants ${\kappa o}_k$ of $\TZ$ where $k \in \Z_8$ \cite{Lin15}.

We list the following presumably difficult problem to understand behaviours of invariants more for Seifert fibered spheres by taking the risk of having negative answers.

\begin{introprob}
For Seifert fibered spheres $Y = \Sigma(a_1, \ldots, a_n)$, what are the possible relations between the following homology cobordism invariants?
\begin{itemize}
\item $\bar{\mu}(Y)$, $w(Y; m)$, and ${\kappa o}_k (Y)$, 
\item $d(Y)$ and $R(Y; e)$.
\end{itemize}
\end{introprob} 

\subsubsection{A Diversification: More about Summands of \texorpdfstring{$\TZ$}{T3Z}}
\label{moresummands}

Around the 2000s, two more epimorphisms of $\TZ$ were found: Ozsv\'ath-Szab\'o $d$-invariant \cite{OS03b} and Fr{\o}yshov $\delta$-invariant\footnote{There are two $h$-invariants of Fr{\o}yshov: the \say{old} one \cite{Fro02} and the \say{new} one \cite{Fro10}. To avoid ambiguity, we follow the notation that appeared in Manolescu's survey \cite{Man20}, called the \say{new} $h$-invariant $\delta$-invariant.} \cite{Fro10}. The latter invariant is also owed to Kronheimer and Mrowka \cite{KM07}. The seminal articles of Kutluhan, Lee, and Taubes \cite{KLT20, KLT20b, KLT20c, KLT20d, KLT20e} yield that $\delta = -d/2$.

Given any relatively coprime positive integers $p,q$ and $r$, the Brieskorn sphere $\Sigma(p,q,r+pq)$ can be obtained by the Brieskorn sphere $\Sigma(p,q,r)$ by applying $(-1)$-surgery along the singular fiber of degree $r$. This topological operation is called \emph{Seifert fiber surgery}, see the paper of Lidman and Tweedy \cite{LT18} for a detailed exposition. 

Performing the above type of Seifert fibered surgeries, the items \eqref{r} and \eqref{s} in Theorem~\ref{summand} can be constructed from the items \eqref{j} and \eqref{k} in Theorem~\ref{summand} respectively. We know that the $d$-invariant remains same under this special Seifert fiber surgery, consult the articles of Lidman and Tweedy \cite{LT18}, Karakurt, Lidman, and Tweedy \cite{KLT21}, and Seetharaman, Yue, and Zhu \cite{SYZ21} for this result. Relying on the computations in \cite{T13} and \cite{KS20} again, we have the following result.
\newpage
\begin{introthm}
\label{summand}
The following homology $3$-spheres individually generate $\Z$ summands in $\TZ$:
\begin{enumerate}
\item \label{j} $\Sigma(p,q,pqn-1)$ for each $n \geq 1$,
\item \label{r} $\Sigma(p,q,+pqn-1+pqm)$ for each $n,m \geq 1$,
\item \label{k} $\Sigma(p_n, q_n, r_n)$ for each $n \geq 1$ where $p_n q_n + p_n r_n - q_n r_n = 1$,
\item \label{s} $\Sigma(p_n, q_n, r_n+p_nq_n m)$ for each $n,m \geq 1$ where $p_n q_n + p_n r_n - q_n r_n = 1$.
\end{enumerate}
\end{introthm} 

In a similar fashion, we can pass to the Brieskorn sphere $\Sigma(p,q,r+2pq)$ from the Brieskorn sphere $\Sigma(p,q,r)$ by applying twice $(-1)$-surgery along the singular fiber of degree $r$. In \cite{SYZ21}, Seetharaman, Yue, and Zhu also observed that the maximal monotone subroots carrying the Floer theoretic invariants do not change after performing the above type of Seifert fiber surgeries consecutively. Recently, in \cite{KS22}, Karakurt and the author presented more families of homology $3$-spheres generating infinite rank summands in $\TZ$ by computing their connected Heegaard Floer homologies \cite{HHL18} effectively and using the invariants of Dai, Hom, Stoffregen, and Truong. Notice that the connected Heegaard Floer homology was introduced by Hendricks, Hom, and Lidman. Further, they proved that it is a homology cobordism invariant itself \cite{HHL18} unlike the classical or involutive Heegaard Floer homology.

Together with the above observation, we can conclude the following theorem. In particular, two collections of families in the items \eqref{l} and \eqref{m} in Theorem~\ref{summands}, and the family of Dai, Hom, Stoffregen, and Truong in Theorem~\ref{DHST} are not homology cobordant to each other for any equal value of $n$, with a single exception, see the discussion in \cite{KS22}. However, their spans in $\TZ$ are not distinct, see \cite[Section~6]{DS17}. 

\begin{introthm}[\cite{DHST18, KS22}]
\label{summands}
The following infinite families of homology $3$-spheres generate $\Z^\infty$ summands in $\TZ$: 
\begin{enumerate}
\item \label{l} $\{\Sigma(2n+1,3n+2,6n+1) \}^\infty_{ n=1}$,
\item \label{m} $\{\Sigma(2n+1,3n+1,6n+5) \}^\infty_{ n=1}$,
\item \label{n} $\{ \Sigma(2n+1,4n+1, 4n+3 + 2m(2n+1)(4n+1)) \}^\infty_{ n,m=1}$,
\item \label{o} $\{ \Sigma(2n+1,3n+2, 6n+1 + 2m(2n+1)(3n+2)) \}^\infty_{ n,m=1}$,
\item \label{p} $\{ \Sigma(2n+1,3n+1, 6n+5 + 2m(2n+1)(3n+1)) \}^\infty_{ n,m=1}$.
\end{enumerate}
\end{introthm} 

\subsection{The Trivial Element of \texorpdfstring{$\TZ$}{T3Z}}
\label{trivial}

A central problem in low-dimensional topology is to investigate the following interaction between $3$- and $4$-manifolds as an algebro-topological analog of the relation between $S^3$ and $B^4$.

\begin{introprob}[Problem~4.2, \cite{K78b}]
\label{bounding}
Which homology $3$-spheres bound contractible $4$-manifolds or homology $4$-balls?
\end{introprob}

There are plenty of examples of Brieskorn spheres that bound Mazur type contractible $4$-manifolds built with a single $0$-, $1$-, and $2$-handle \cite{M61}. Following Kirby's celebrated work \cite{K78}, some classical articles appeared subsequently: Akbulut and Kirby \cite{AK79}, Casson and Harer \cite{CH81}, Stern \cite{S78}, Fintushel and Stern \cite{FS81}, Maruyama \cite{M81, M82}, and Fickle \cite{F84}. In addition, some of these results were found independently of Kirby calculus, see Fukuhara \cite{F78} and Martin \cite{M79}. Some of these families also bound Po\'enaru manifolds, contractible $4$-manifolds built with a $0$-handle, many $1$- and $2$-handles, see \cite{P60, S20b, AS22}.

\begin{introthm}
The following homology $3$-spheres bound Mazur manifolds with one $0$-handle, one $1$-handle and one $2$-handle. Further, $\Sigma(2,7,47)$ and $\Sigma(3,5,49)$ bound homology $4$-balls.
\begin{itemize}
\item $\Sigma(2,3,13)$, $\Sigma(2,3,25), \Sigma(2,7,19)$, $\Sigma(3,5,19)$.
\item $\Sigma(p,ps-1,ps+1)$ for $p$ even and $s$ odd,
\item $\Sigma(p,ps \pm 1,ps \pm 2)$ for $p$ odd and $s$ arbitrary,
\item $\Sigma(2,2s \pm 1, 2 \cdot 2 \cdot (2s \pm 1) + 2s \mp 1)$ for $s$ odd, 
\item $\Sigma(3,3s \pm 1, 2 \cdot 3 \cdot (3s \pm 1) + 3s \mp 2)$ for $s$ arbitrary, 
\item $\Sigma(3,3s \pm 2, 2 \cdot 3 \cdot (3s \pm 2) + 3s \mp 1)$ for $s$ arbitrary.
\end{itemize}
\end{introthm} 

It would be interesting to compare the existence of homology $3$-spheres bounding contractible $4$-manifolds and homology $4$-balls, so we may address the following problem. The possible candidates for Seifert fibered spheres are two examples of Fickle: $\Sigma(2,7,47)$ and $\Sigma(3,5,49)$. They are known to bound only homology $4$-balls.

\begin{introprob}
\label{compare}
Is there any Seifert fibered sphere $\Sigma(a_1, \ldots, a_n)$ which bounds a homology $4$-ball but not a contractible $4$-manifold?
\end{introprob}

Note that Problem~\ref{compare} is known for $\Sigma(2,3,5) \# -\Sigma(2,3,5)$.\footnote{In general, it is known for a homology $3$-sphere which bounds a simply-connected $4$-manifold with non-standard definite intersection form. Taubes attributed this result to Akbulut.} It cannot bound a contractible $4$-manifold, see Taubes' article \cite[Proposition~1.7]{T87}. However, the isomorphism of Gonz\'{a}lez-Acu\~{n}a in Theorem~\ref{Gonz} guarantees that every homology $3$-sphere bounding homology $n$-ball automatically bounds contractible $n$-manifold unless $n=3$.

When the number of fibers increases, there is a bold conjecture, which was first indicated by Fintushel-Stern, explicitly stated by Lawson \cite{L88}, and later highlighted by Koll\'ar \cite[Conjecture~20]{K08}. This problem is closely related to the Montgomery-Yang problem motivated by the previous results in both algebraic geometry and gauge theory. The problem expects that every pseudo-free circle action on the $5$-dimensional sphere has at most $3$ non-free orbits \cite[Conjecture~6]{K08}. Note that some computational verifications of this conjecture was provided in the paper of Lawson \cite{L88}.

\begin{introprob}[Three Fibers Conjecture]
\label{TFC}
Is there any Seifert fibered sphere $\Sigma(a_1, \ldots, a_n)$ with $n>3$ which bounds a homology $4$-ball?
\end{introprob}

Problem~\ref{TFC} cannot be generalized for plumbed homology $3$-spheres that are not Seifert fibered.\footnote{\label{xxx} Note that $\partial X(1) = \Sigma(2,5,7)$ and $\partial X'(1) = \Sigma(3,4,5)$, compare with \cite{AK79}, \cite{CH81}, and \cite{S20b}. Therefore, they are not Seifert fibered unless $n=1$.} The first examples were given by Maruyama \cite{M82}, independently obtained by Akbulut and Karakurt \cite[Theorem~1.4]{AK14}. In \cite{S20b}, we presented two more family of plumbed homology $3$-spheres bounding contractible $4$-manifolds. 

\begin{figure}[ht]
\begin{center}
\includegraphics[width=.6\textwidth]{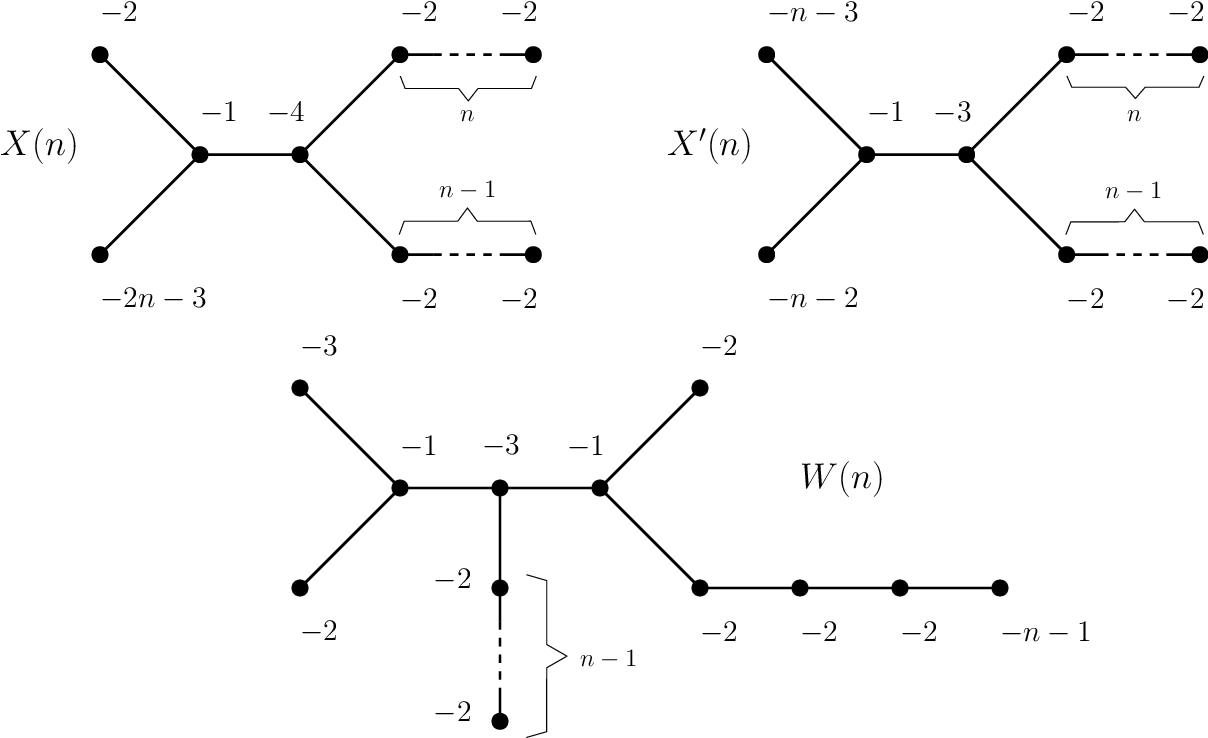}   
\end{center}
\caption{The plumbing graphs of $X(n)$, $X'(n)$, and $W(n)$.}
\label{fig:plumbed}
\end{figure}

\begin{introthm}[Theorem~1, \cite{M82}; Theorem~1.4-5, \cite{S20b}]
Let $X(n)$, $X'(n)$, and $W(n)$ be Maruyama, the companion of Maruyama, and Ramanujam plumbed $4$-manifold, shown in Figure~\ref{fig:plumbed}. Then for each $n\geq 1$, boundaries $\partial X(n)$ and $\partial X'(n)$ bound  Mazur manifolds with one $0$-handle, one $1$-handle and one $2$-handle. Further, the boundary of $\partial W(n)$ bounds a Po\'enaru manifold with one $0$-handle, two $1$-handles and two $2$-handles for $n\geq 1$.
\end{introthm}

Note that $W(1)$ is known as \emph{Ramanujam surface}, the famous \emph{homology plane} constructed by Ramanujam \cite{R71}. It is the first example of an algebraic complex smooth surface sharing the same of homology of the complex plane $\mathbb{C}^2$ but not analytically isomorphic to $\mathbb{C}^2$. We call a non-trivial homology $3$-sphere a \emph{Kirby-Ramanujam spheres} if it bounds both a homology plane and a Mazur/Po\'enaru type contractible $4$-manifold. In \cite{AS22}, Aguilar and the author found several infinite families of Kirby-Ramanujam spheres in the light of Problem~\ref{bounding}.

In \cite{A91}, Akbulut introduced very crucial geometric objects called \emph{corks}. These are defined to be contractible smooth $4$-manifolds together with involutions on the boundary $3$-manifolds, which extend to self-homeomorphisms but not to self-diffeomorphisms of the ambient manifolds. As they generate all exotic phenomena for simply-connected $4$-manifolds via cork twists \cite{CFHS96, M96}, they draw a special interest in low-dimensional topology. Corks have recently been studied extensively using Heegaard Floer homology by Dai, Hedden and Mallick \cite{DHM20}, and they introduced an algebraic object called \emph{homology bordism group of involutions} $\Theta^\tau_\Z$ as a modification of the homology cobordism group $\TZ$. However, the following question remains a very interesting open problem:

\begin{introprob}
Is there any Seifert fibered space $\Sigma(a_1, \ldots, a_n)$ bounding a cork? 
\end{introprob}

Seifert fibered spaces cannot appear as the boundaries homology planes due to Orevkov \cite{O97}. However, the splice of Seifert manifolds along their singular fibers are shown to bound homology planes \cite{AS22}. Since they also bound contractible $4$-manifolds, we can pose the following problem. If such a homology $3$-sphere exists, then after possibly applying cork twists, we can glue these contractible $4$-manifolds along their common boundary. This gives a homotopy $4$-sphere so that it is homeomorphic to the $4$-sphere $S^4$ by Freedman \cite{F82}. Therefore, this $4$-manifold would be a new potential candidate counterexample to the smooth Poincar\'e Conjecture in dimension $4$. 

\begin{introprob}
Is there any homology $3$-sphere bounding both a cork and a contractible homology plane? 
\end{introprob}

Using the surgery descriptions of $\Sigma(p,q,pq \mp 1)$ in terms of torus knots, one can prove the following theorem as an immediate corollary of the main results of Gordon \cite{G75} and Lidman, Karakurt, and Tweedy \cite{KLT21}. For the constructive part, an alternative direct proof can be given by finding the plumbing graphs of splices explicitly \cite{EN85} and doing Kirby calculus. The obstruction of knots bounding smooth disks requires the result of Lidman and Tweedy \cite{LT18}.

\begin{introthm}
Let $K(pq \mp 1)$ denote the singular fiber in $\Sigma(p,q,pq \mp 1 )$. Then $K(pq \mp 1)$ is not smoothly slice in $\Sigma(p,q,pq \mp 1)$, and $\Sigma(p,q,pq \mp 1 )$ does not bound a contractible $4$-manifold. However, the following splicing homology $3$-spheres bound Po\'enaru manifolds with one $0$-handle, $p$ $1$-handles, and $p$ $2$-handles: \begin{center} $\Sigma(p,q,pq-1) \tensor[_{K(pq-1)}]{\bowtie}{_{K(pq+1)}} \Sigma(p,q,pq+1).$ \end{center}
\end{introthm}

Independent results of Hirsch, Rokhlin, and Wall around the 1960s indicate that every homology $3$-sphere is smoothly embedded in $S^5$, see \cite{H61}, \cite{R65} and \cite{W65}. Making the target space smaller,  we may ask which homology $3$-spheres can be embedded in $S^4$. In the topological category, the problem has a complete answer thanks to Freedman's celebrated article \cite{F82}: every homology $3$-sphere is topologically embedded in $S^4$. Adding an extra smoothness condition, we can state another wide open problem in low-dimensional topology.

\begin{introprob}[Problem~3.20, \cite{K78b}]
\label{S4}
Which homology $3$-spheres can be smoothly embedded in $S^4$?
\end{introprob}

Another simple algebraic topology observation indicates that a homology $3$-sphere smoothly embedded in $S^4$ splits $S^4$ into two homology $4$-balls. Therefore, homology cobordism invariants provide obstructions for the smooth embeddings of homology $3$-spheres in $S^4$.

One can wonder about the reverse direction of the above observation. Studying branched coverings of cross sectional slices of knotted $2$-spheres $S^2$ in $S^4$, McDonald provided the first examples of homology $3$-spheres which are smoothly embedded in a homology $4$-ball but not any homotopy $4$-sphere \cite{Mc22}. His examples are certain double cyclic branched coverings of spuns of torus knots. We may address this implication to Seifert fibered manifolds and ask:

\begin{introprob}
Is there any Seifert fibered sphere which bounds a homology $4$-ball but cannot be smoothly embedded in $S^4$?
\end{introprob}

\subsection{Generators of \texorpdfstring{$\TZ$}{T3Z}}
\label{generator}

The first result concerning the generators of $\TZ$ was owed to Freedman and Taylor.

\begin{introthm}[Corollary~1B, \cite{FT77}]
The group $\TZ$ is generated by homology $3$-spheres which are boundaries of $4$-manifolds having the homology of $S^2 \times S^2$.
\end{introthm}

A homology $3$-sphere $Y$ is called \emph{irreducible}\footnote{A homology $3$-sphere $Y$ is said to be \emph{prime} if it cannot be written as a connected sum of two homology $3$-spheres non-trivially (i.e. either summand is not $S^3$). For homology $3$-spheres, sometimes the terms of \enquote{prime} and \enquote{irreducible} can be used interchangeably unless $Y=S^3$, see \cite[Lemma~1]{M62}.} if every embedded $2$-sphere $S^2$ in $Y$ is the boundary of an embedded $B^3$. Livingston showed that irreducible homology $3$-spheres are generic enough to generate the homology cobordism group.

\begin{introthm}[Theorem~3.2, \cite{L81}]
\label{L}
Every class in $\TZ$ admits an irreducible representative.
\end{introthm}

We call a homology $3$-sphere $Y$ \emph{hyperbolic} if $Y$ is a geodesically complete Riemannian $3$-manifold of constant sectional curvature $-1$. The geodesically completeness requires that at any point $p \in Y$, the geodesic exponential map $\mathrm{exp}_p$ on $T_pY$ is the entire tangent space at $p$. Myers proved that every homology cobordism class admits a hyperbolic representative.

\begin{introthm}[Theorem~5.1, \cite{M83}]
\label{My}
Every class in $\TZ$ admits a hyperbolic representative.
\end{introthm}

A pair $(Y,\xi)$ is called \emph{Stein fillable} if there is a Stein domain $(X,J,\phi)$ where $\phi$ is bounded below, $Y$ is an inverse image of an regular value of $\phi$, and $\xi = \mathrm{ker}(-d\phi \circ J)$ is an induced contact structure. Mukherjee showed that the generator set of $\TZ$ can be chosen as Stein fillable homology $3$-spheres \cite{M20}.

\begin{introthm}[Theorem~1.5, \cite{M20}]
The group $\TZ$ is generated by Stein fillable homology $3$-spheres.
\end{introthm}

In contrast to the above positive directional results, various computations of homology cobordism invariants of homology $3$-spheres lead to the following observation of Fr{\o}yshov \cite{F16}, Stoffregen \cite{St17}, Lin \cite{L17}, and Nozaki, Sato, and Taniguchi \cite{NST22}.

\begin{introthm}
There exist several infinite families of homology $3$-spheres that are not homology cobordant to any Seifert fibered homology sphere.
\end{introthm}

In \cite{HHSZ20}, Hendricks, Hom, Stoffregen, and Zemke established a surgery exact triangle formula for the involutive Heegaard Floer homology. As an application, they provided a homology $3$-sphere not homology cobordant to any linear combination of Seifert fibered spheres, \cite[Theorem~1.1]{HHSZ20}. This manifold is obtained by integral Dehn surgery on a combination of torus knots and a cable of a torus knot: $S^3_{+1}(-T_{6,7} \# T_{6,13} \# -T_{2,3;2,5})$. Hence, Seifert fibered manifolds are not generic enough to generate $\TZ$:

\begin{introthm}[Theorem~1.1, \cite{HHSZ20}]
\label{HHSZ}
The Seifert fibered spheres cannot generate the group $\TZ$. Therefore, $\Theta^3_{SF}$ is a proper subgroup of $\TZ$. Further, $\TZ / \Theta^3_{SF}$ has a $\Z$ subgroup.
\end{introthm}

Here, $\Theta^3_{SF}$ denotes the subgroup of $\TZ$ generated by Seifert fibered spheres. Note that $S^3 = \Sigma(1,q,r)$. By using Kirby calculus, Nozaki, Sato, and Taniguchi proved that the example of Hendricks, Hom, Stoffregen, and Zemke is a graph homology $3$-sphere, see \cite[Appendix~A]{NST22}. Therefore, we can ask the following question as to the next step of obstructions:

\begin{introprob}
Do graph homology $3$-spheres generate the group $\TZ$?
\end{introprob}

Let $\Theta^3_{G}$ denote the subgroup of $\TZ$ generated by graph homology $3$-spheres. The previous problem is equivalent to asking whether $\Theta^3_{G} = \TZ$ or not. Nozaki, Sato, and Taniguchi proposed a strategy in \cite[Conjecture~1.19]{NST22} so that likely $\Theta^3_{G} \lneq \TZ$.

Hendricks, Hom, Stoffregen, and Zemke compared the subgroup $\Theta^3_{SF}$ with the whole group $\TZ$ in another work and they were able to provide the existence of an infinitely generated subgroup in the quotient $\TZ / \Theta^3_{SF}$ spanned by the family of homology $3$-spheres $S^3_{+1}(-T_{2,3} \# -2T_{2n,2n+1} \# -T_{2n,4n+1})$ for odd $n \geq 3$:

\begin{introthm}[Theorem~1.1, \cite{HHSZ21}]
\label{HHSZ2}
The quotient $\TZ / \Theta^3_{SF}$ has a $\Z^\infty$ subgroup.
\end{introthm}

The new immediate challenge would be to ask:

\begin{introprob}
Does the quotient $\TZ / \Theta^3_{SF}$ contain a $\Z^\infty$ summand?
\end{introprob}

Another curiosity about the possible generators of $\TZ$ is of course surgeries on knots in the $3$-sphere. One can expect that these manifolds are not sufficient to provide a generating set for $\TZ$, see \cite[Corollary~1.7]{NST22}. However, the following problem still remains open.

\begin{introprob}
Do surgeries on knots in $S^3$ generate $\TZ$?
\end{introprob}

\subsection{Torsion of \texorpdfstring{$\TZ$}{T3Z}} 
\label{torsion}

In their seminal articles, Matumoto \cite{M78} and Galewski and Stern \cite{GS80} reduced the triangulation conjecture to a problem about the interplay between $3$- and $4$-manifolds up to homology cobordism. Since then $\TZ$ has been a very attractive object in low-dimensional topology.
A splitting would provide a homology $3$-sphere $Y$ such that $\mu(Y) = 1$ and $Y$ is $2$-torsion in the homology cobordism group.

\begin{introthm}[\cite{M78, GS80}]
\label{GSM}
For $n \geq 5$, there exist non-triangulable topological $n$-manifolds if and only if the following exact sequence does not split \begin{equation} \label{split} 0 \longrightarrow \mathrm{ker}(\mu) \longrightarrow \TZ \overset{\mu} \longrightarrow \Z_2 \longrightarrow 0. \tag{$\star$}
\end{equation}
\end{introthm}

Prior to the work \cite{M78, GS80}, Casson asked whether every homology $3$-sphere $Y$ with an orientation reversing diffeomorphism satisfies $\mu(Y) = 0$, see \cite[Problem~3.43]{K78b}. If it were false, then $Y \# Y = Y \# -Y$ would bound the homology $4$-ball $(Y \setminus \intE) \times [0,1]$, giving an element of order $2$ in $\TZ$. Independently, Birman (in an unpublished note), Galewski and Stern \cite{GS79}, and Hsiang and Pao \cite{HP79} partially answered this question affirmatively for homology $3$-spheres with orientation-reversing involutions. Finally, Casson showed that the $\mu$-invariant must be zero for such a homology $3$-sphere $Y$ in general \cite{AM90}.

Next, Saveliev \cite{Sav02c} proved that $\Z_2$ torsion in the homology cobordism group cannot be generated by Seifert fibered spaces (plumbing homology $3$-spheres in general) with non-trivial Rokhlin invariants. He showed that such a Seifert manifold must be of infinite order by extending the previous work of Fukumoto, Furuta, and Ue \cite{FFU01}.

Finally, Manolescu \cite{M16} constructed $\mathrm{Pin} (2)$-equivariant Seiberg-Witten Floer homology and provided three sensitive invariants of homology $3$-spheres. They are called $\alpha, \beta,$ and $\gamma$ invariants of $\TZ$. Specifically, the Manolescu $\beta$-invariant has the following three crucial properties:
\begin{enumerate}
\item \label{minus} $\beta(-Y) = -\beta(Y)$,
\item \label{lift} $- \beta(Y) = \mu(Y) \mod 2$ where $\mu$ is the Rokhlin invariant,
\item \label{Manolescu} $\beta$ is an invariant of $\TZ$.
\end{enumerate}

The existence of the Manolescu $\beta$-invariant guaranteed that the exact sequence \eqref{split} does not split and leads to the disproof of the triangulation conjecture, see \cite[Problem~4.4]{K78b} and \cite{Man16, M16, M18}. For this achievement, the homology cobordism invariance of the Manolescu $\beta$-invariant is particularly critical because beforehand there exist invariants satisfying properties both \eqref{minus} and \eqref{lift} but not \eqref{Manolescu}; for instance, the Casson invariant $\lambda$. Therefore, it cannot be used for the rejection of the triangulation conjecture for high-dimensional manifolds; however, it is sufficient for disproval of the conjecture for the particular case of $n=4$. See the book of Akbulut and McCarthy \cite{AM90} for details. For an alternative disproof of the triangulation conjecture for high-dimensional manifolds using a similar strategy, see F. Lin's monograph \cite{L18}.

Since the Manolescu $\beta$ invariant provides an integral lift of the Rokhlin invariant $\mu$, he also ruled out the existence of $\Z_2$ torsion in $\TZ$ for the following type of homology $3$-spheres: 

\begin{introthm}[Corollary~1.2, \cite{M16}]
\label{M}
Let $Y$ be a homology $3$-sphere such that $\mu(Y)=1$. Then $Y$ cannot represent $\Z_2$ torsion in $\TZ$. In other words, $Y \# Y$ cannot bound a homology $4$-ball.
\end{introthm}

Currently, we do not know whether there exists a non-trivial homology $3$-sphere $Y$ with a vanishing $\mu$-invariant  so that $Y \# Y$ bounds a contractible $4$-manifold or a homology $4$-ball. Also, we have no further obstructions for other types of torsion in $\TZ$. Hence we curiously state the following problem:

\begin{introprob}
\label{torprob}
Does the group $\TZ$ contain any torsion $\Z_n$ for $n \geq 2$? Modulo torsion, is $\TZ$ free abelian?
\end{introprob}

Only for the $\Z_2$ type torsion, there are some new candidates found in the recent work of Boyle and Chen \cite{BC22}. These examples originate from cyclic double branched coverings of $S^3$ along certain non-slice strongly negative amphichiral knots of determinant $1$.

\section{Two Relatives of \texorpdfstring{$\TZ$}{T3Z}}

Finally, we discuss the close and crucial relationship between the knot concordance group $\mathcal{C}$, the homology cobordism group $\TZ$, and the rational homology cobordism group $\TQ$.

\subsection{The Elder: The Knot Concordance Group \texorpdfstring{$\mathcal{C}$}{C}}
\label{Cconcordance}

A \emph{knot} $K$ is a smooth embedding of a circle $S^1$ into $S^3$. The \emph{knot concordance group} $\mathcal{C}$ is defined as \begin{center}
$\mathcal{C} = \{ \text{oriented knots up to isotopy} \} / \sim$ \end{center} where the equivalence relation \emph{concordance} $\sim$ is given for two arbitrary knots $K_0$ and $K_1$ as 

\vspace{0.5 em}

\noindent \begin{minipage}{0.7\linewidth}

\[  K_0 \sim K_1 
 \iff  \begin{cases}
\ \text{There exists a cylinder} \ C \ \text{such that} \\ 
   \bullet \  C \subset S^3 \times [0,1], \\
   \bullet \ \partial C = -(K_0) \cup K_1.
  \end{cases} \]
  
\end{minipage}  
\begin{minipage}{0.3\linewidth}  
  
\includegraphics[width=0.87\textwidth]{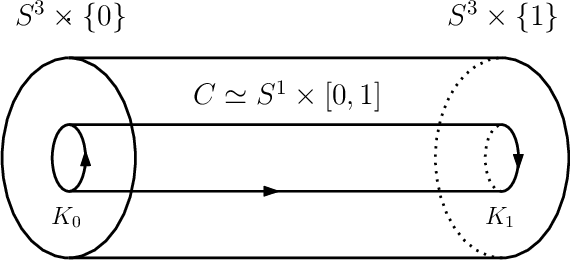}   
\end{minipage}  

\vspace{0.8 em}

Fox and Milnor introduced the group $\mathcal{C}$ in their celebrated article \cite{FM66}. The summation is induced by connected sum of knots. The concordance class of the unknot gives the zero element. Inverse elements are found by mirroring knots and reversing their orientations.

Knots concordant with the unknot are said to be \emph{slice knots}. Equivalently, slice knots are the knots that bound smoothly embedded disks in $B^4$. \emph{Ribbon knots} can be defined by restricting of the handle decomposition of the smooth disks; they are the ones that bound such disks without $2$-handles. Clearly, every ribbon knot is a slice. However, the opposite is one of the most famous long-standing problems in knot theory proposed by Fox \cite{F62}:

\begin{introprob}[Slice-Ribbon Conjecture]
Is every slice knot is ribbon?
\end{introprob}

There are candidates for a counterexample to the slice-ribbon conjecture, provided by Gompf, Scharlemann and Thompson \cite{GST10} and Abe and Tagami \cite{AT16}. On the other hand, this conjecture was confirmed for $2$-bridge knots by Lisca \cite{L07, L07b} and for most pretzel and Montesinos knots by Greene and Jabuka \cite{GJ11} and Lecuona \cite{Lec12, Lec15, Lec18, Lec19}.

In his celebrated work \cite{G81}, Gordon defined the notion of \emph{ribbon concordance} as an analogue of ribbon knots so that the Morse function induced by the concordance $S^3 \times [0,1] \to [0,1]$ has no critical points of index 2. Furthermore, Gordon conjectured that the ribbon concordance is a partial order; this was recently proved by Agol \cite{A22}. Zemke \cite{Z19} initiated an approach to the study of ribbon concordance using knot Floer homology, which was generalized to $3$-manifolds by Daemi, Lidman, Vela-Vick, and Wong \cite{DLVW22}. Their formalism also provides important links to Thurston geometries. 

A careful analysis of the classical articles of Fox and Milnor \cite{F62, FM66}, Murasugi \cite{M65}, Robertello \cite{Rob65}, Levine \cite{L69} and Tristam \cite{T69} ensured the existence of infinitely generated $\Z^\infty$ and $\Z_2^\infty$ summands of the knot concordance group so that we pose the following first question regarding the algebraic structure of $\mathcal{C}$: 

\begin{introprob}
Is the group $\mathcal{C}$ isomorphic to $\Z^\infty \oplus \Z_2^\infty$?
\end{introprob}

Levine's eminent articles provide a surjective homomorphism $\phi: \mathcal{C} \to \Z^\infty \oplus \Z_2^\infty \oplus \Z_4^\infty$ \cite{L69, L69b}. First, Casson and Gordon \cite{CG78} proved that $\phi$ is not an isomorphism. Next, Jiang \cite{J81} improved their result by showing that $\mathrm{Ker}(\phi)$ has a $\Z^\infty$ subgroup. Finally, Livingston exhibited that $\mathrm{Ker}(\phi)$ has a $\Z_2^\infty$ subgroup \cite{L01}. The following question remains open:

\begin{introprob}
\label{Levine}
Does Levine's homomorphism $\phi$ split?
\end{introprob}

An affirmative answer to Problem~\ref{Levine} will provide elements of order $4$ in $\mathcal{C}$. Furthermore, it will guarantee that elements of order $2$ do not arise only from negative amphicheiral knots, see \cite{L05} for more details. Furthermore, obstructions to elements of order $4$ were found by Livingston and Naik \cite{LN99}. Therefore, Problem ~\ref{Levine} is closely related to the remaining finite part of the knot concordance group.

\begin{introprob}
Does the group $\TZ$ contain any torsion $\Z_n$ for $n > 2$?
\end{introprob}

In \cite{COT03, COT04}, Cochran, Orr and Teichner introduced and studied the deep structure of $\mathcal{C}$ by forming a filtration of the group via an infinite sequence of subgroups \begin{center}
$ \ldots \subset \mathcal{F}_{n+1} \subset \mathcal{F}_{n.5} \subset \mathcal{F}_{n} \subset  \ldots  \subset \mathcal{F}_{1.5} \subset \mathcal{F}_{1} \subset \mathcal{F}_{0.5} \subset \mathcal{F}_{0} \subset \mathcal{C} $
\end{center}

\noindent where $\mathcal{F}_{0}$, $\mathcal{F}_{0.5}$, and $\mathcal{F}_{1.5}$ respectively correspond to knots with trivial Arf invariant, knots in the kernel of $\phi$, and knots having vanishing Casson–Gordon invariants. This filtration structure is highly non-trivial; in particular, Cochran, Harvey, and Leidy proved that each quotient $\mathcal{F}_{n} / \mathcal{F}_{n.5} $ contains a $\Z^\infty \oplus \Z_2^\infty$ subgroup \cite{CHL09, CHL11}.

The group $\mathcal{C}$ and $\TZ$ are related by the following maps \begin{center} $S^3_{1/n}: \mathcal{C} \to \TZ, \ \ \ [K] \mapsto [S^3_{1/n}(K)].$\end{center} These maps are not homomorphisms but send identity to identity, see classical work of Gonz\'{a}lez-Acu\~{n}a \cite{G70} and Gordon \cite{G75}.

The set of maps $S^3_{1/n}$ was used by Peters to study the knot concordance with the help of the Heegaard Floer theoretic $d$-invariant \cite{P10}. The same technique was adapted in the work of Hendricks and Manolescu \cite{HM17} in the setup of involutive Heegaard Floer homology. This approach can be applied a priori to the other homology cobordism invariants.

Finally, we briefly mention key obstructive techniques originating from several theories of knots, $3$- and $4$-manifolds. Akbulut and Matveyev \cite{AM97} and Rudolph \cite{R95} used contact geometry in the spirit of Eliashberg's work \cite{E90}. The gauge theoretic methods of Donalson and Taubes \cite{D83, T87} were adapted by Cochran and Gompf \cite{CG86}, Fintushel and Stern \cite{FS85}. Casson-Gordon invariants \cite{CG78, CG86} were applied successfully by Litherland \cite{L84}, Kirk and Livingston \cite{KiL99}, Friedl \cite{F04},  Kim\cite{K05}, and Aceto, Golla, and Lecuona \cite{AML18}. The knot Floer homology independently defined Ozsv\'ath and Szab\'o  \cite{OS04b} and Rasmussen \cite{R03} has been used extensively, see for example Ozsv\'ath and Szab\'o \cite{OS03f} and Ozsv\'ath, Szab\'o, and Stipsicz \cite{OSS17}. Furthermore, Khovanov homology and Lee's refinement \cite{K00, L05} provided powerful invariants and techniques through the work of Rasmussen \cite{R10}, Kronheimer and Mrowka  \cite{KM13}, Lipshitz and Sarkar \cite{LS14}, and Piccirillo \cite{P20}. Recently, Dai, Hom, Stoffregen, and Truong produced involutive Floer theoretic invariants \cite{DHST19}, building on the work of Hendricks and Manolescu \cite{HM17}. Moreover, Khovanov-Rozansky homology \cite{KR08} was used by Lobb \cite{L09} and Lewark \cite{L14} to provide quantum obstructions. Finally, instanton knot Floer homology \cite{Fl90} has yielded crucial results led by Kronheimer and Mrowka \cite{KM10, KM11}, Hedden and Kirk \cite{HK12}, and Baldwin and Sivek \cite{BS21, BS22}.

For more details and further advancements, see survey of Gordon \cite{G78}, Livingston \cite{L05}, Hom \cite{H17, H21}, and problem lists \cite{Pr11, Pr16, Pr19}.

\subsection{The Younger: The Rational Homology Cobordism Group \texorpdfstring{$\TQ$}{T3Q}}
\label{Qcobordism}

Changing the role of integer coefficients with rational ones in the definition of $\TZ$, we obtain the \emph{rational homology cobordism group} $\TQ$. Deciphering the trivial class of this group has been of special interest in low-dimensional topology, constituting the following problem attributed to Casson:

\begin{introprob}[Problem~4.5, \cite{K78b}]
\label{bounding2}
Which rational homology $3$-spheres bound rational homology $4$-balls?
\end{introprob}

From both constructive and obstructive perspectives, Problem~\ref{bounding2} has been studied extensively with the help of the techniques introduced by Casson and Gordon \cite{CG78}. For each prime $p$ and $r\geq 1$, we have a group homomorphism\begin{center}
 $\Sigma_{p^r}: \mathcal{C} \to \TQ , \ \ \ [K] \mapsto [\Sigma_{p^r}(K)].$
\end{center}

The homomorphism of Casson and Gordon was used for the construction of concordance invariants. See the work of Manolescu and Owens \cite{MO07}, Jabuka \cite{J12}, Alfieri, Kang and Stipsicz \cite{AKS19}, and Baraglia \cite{Bar22}.

The work of Lisca \cite{L07, L07b} on the slice-ribbon conjecture for $2$-bridge knots led to the classification of lens spaces and sums of lens spaces bounding rational homology $4$-balls. Similarly, the articles of Greene and Jabuka \cite{GJ11} and Lecuona \cite{Lec12, Lec15, Lec18, Lec19} provided Seifert fibered rational homology $3$-spheres bounding rational homology $4$-balls. Recently, Aceto and Golla \cite{AG17} and Aceto, Golla, Larson, and Lecuona \cite{AGLL20} classified surgeries on torus knots that bound rational balls. Also, Lokteva \cite{Lis20} extended their results to cables of torus knots. Furthermore, Maruyama \cite{Mar80}, Fintushel and Stern \cite{FS80}, Casson and Harer \cite{CH81}, Etnyre and Tosun \cite{ET20}, Simone \cite{Sim20, Sim20b} and Lokteva \cite{Lis22} constructed various rational homology $3$-spheres bounding rational homology $4$-balls by using Kirby calculus and knot theory; see also \cite{L07, L07b, Lec12, AGLL20} for the construction of certain spaces.
 
Several theories extended to rational homology $3$-spheres and their invariants can be extensively used for powerful obstructions. Consult the articles by Owens and Strle \cite{OStr06}, Simone \cite{Sim20b}, Choe and Park \cite{CP21}, and Greene and Owens \cite{GO22} using Donaldson's diagonalization theorem and Heegaard Floer homology; Casson and Gordon \cite{CG86}, Fintushel and Stern \cite{FS87b}, Mati\'c \cite{Mat88}, Ruberman \cite{R88}, Yu \cite{Y91}, and Mukawa \cite{M02} using Casson-Gordon invariants and gauge theory; Wahl \cite{W81, W11}, Stipsicz, Szab\'o, and Wahl \cite{SSW08}, and Bhupal and Stipsicz \cite{BS11} using singularity theory; Baraglia and Hekmati using Seiberg-Witten-Floer theory \cite{BH21, BH22}.  

A combination of the classical work of Casson and Harer \cite{CH81} and Litherland \cite{L79} indicate that $\mathrm{Ker}(\Sigma_{p})$ contains a $\Z^\infty$ subgroup for any prime $p$. In particular, Aceto and Larson showed that $\mathrm{Ker}(\Sigma_{2})$ has a $\Z^\infty$ summand. Further, Aceto, Celoria, and Park \cite{ACP20} proved that $\mathrm{Coker}(\Sigma_{p^r})$ contains a subgroup isomorphic to $\Z^\infty$ if $p \equiv 3 (\mathrm{mod} 4)$ and $\Z^\infty \oplus \Z_2^\infty$ otherwise.

\begin{introprob}
Describe other types of subgroups or summands of $\mathrm{Ker}(\Sigma_{p^r})$ and $\mathrm{Coker}(\Sigma_{p^r})$.
\end{introprob}

In particular, the linear independence of collections of rational homology $3$-spheres in $\TQ$ has been studied by Hedden, Livingston, and Ruberman \cite{HLR12} and Golla and Larson \cite{GL21} using Heegaard Floer homology. See also the work of Mukawa \cite{M02} in the machinery of gauge theory. Nevertheless, the detection of summands in the rational homology cobordism group is an open problem.

\begin{introprob}
Does the group $\TQ$ contain a $\Q^n$ summand for $n \geq 1$? 
\end{introprob}

When Lisca classified connected sums of lens spaces bounding rational homology $4$-balls \cite{L07b}, and he found $2$-torsion elements in $\TQ$. However, the existence of other types of torsion is currently unknown.

\begin{introprob}
Does the group $\TQ$ contain any $n$-torsion for $n > 2$? 
\end{introprob}

We have a natural group homomorphism \begin{center} $\psi: \TZ \to \TQ$
\end{center} induced by inclusion. It is known that the map $\psi$ is not injective. There exists homology $3$-spheres listed in Theorem~\ref{ratballs} that represent non-trivial elements in $\mathrm{Ker}(\psi)$ by the work of Fintushel and Stern \cite{FS84}, Akbulut and Larson \cite{AL18}, the author \cite{S20}, and Simone \cite{Sim20}.\footnote{Note that these families of Brieskorn spheres all bound rational homology $4$-balls for all values of $n$. Simone's family can be generalized in the sense that $S^3_{-1} (K)$ (resp. $S^3_{+1} (K)$) bounds a rational homology $4$-ball when $K$ is an unknotting number one knot with a positive (resp. negative) crossing that can be switched to unknot $K$.} 

\begin{introthm}
\label{ratballs}
The following homology $3$-spheres bound rational homology $4$-balls but do not bound homology $4$-balls. Therefore, they non-trivially lie in $\mathrm{Ker}(\psi)$ since  they all have non-vanishing Rokhlin invariant:
\begin{enumerate}
\item $\Sigma(2,3,7)$, $\Sigma(2,3,19)$,
\item $\Sigma(2,4n+1,12n+5)$, $\Sigma(3,3n+1,12n+5)$ for odd $n \geq 1$,
\item $\Sigma(2,4n+3,12n+7)$, $\Sigma(3,3n+2,12n+7)$ for even $n \geq 2$,
\item \label{simone} $S^3_{-1}(K_n)$ where $K_n$ is the twist knot for odd $n \geq 1$.
\end{enumerate}
\end{introthm}

Furthermore, $\mathrm{Ker}(\psi)$ has a $\Z$ subgroup generated by any single homology $3$-sphere listed above except those in \eqref{simone} because they have non-zero $\bar{\mu}$-invariants. In particular, $\mu$-invariants of Simone's examples in the item~\eqref{simone} are non-trivial. One can expect that $\mathrm{Ker}(\psi)$ might be larger than $\Z$, including some linearly independent infinite subset of these homology $3$-spheres. Thus, we ask the following problem, first posed by Akbulut and Larson \cite{AL18}:

\begin{introprob}
Does $\mathrm{Ker}(\psi)$ contain $\Z^\infty$ subgroup or $\Z^\infty$ summand?
\end{introprob}

It is worthwhile to note that all current homology cobordism invariants cannot detect the linear independence of Brieskorn spheres listed in Theorem~\ref{ratballs} in $\TZ$, see the discussion in Subsection~\ref{recovery}. This is also true for Simone's family, see surgery formulae of the relevant homology cobordism invariants.

The existence of these homology $3$-spheres has a nice application in symplectic geometry. Let $(X,\omega)$ be a symplectic $4$-manifold. A \emph{Stein domain} is a triple $(X,J,\phi)$ such that $J$ is complex structure on $X$ and $\phi: X \to \mathbb{R}$ is a proper plurisubharmonic function. Here, $\phi$ provides compact level sets and a symplectic form: $\phi$ is smooth such that $\phi^{-1}((-\infty,c])$ is compact for all $c \in \mathbb{R}$ and $\omega_\phi (v,w) = -d(d\phi \circ J)(v,w)$ gives a symplectic form. The handle decompositions of Stein domains are completely characterized in the celebrated articles of Eliashberg \cite{E90} and Gompf \cite{G98}: A $4$-manifold is a Stein domain if and only if it has a handle decomposition with $0$-handles, $1$-handles, and $2$-handles; and the $2$-handles are attached along Legendrian knots with framing $\mathrm{tb} - 1$, where $\mathrm{tb}$ denotes the Thurston-Bennequin number. 

If we choose any homology $3$-sphere listed in Theorem~\ref{ratballs}, then the handle decomposition of the corresponding rational ball must contain $3$-handles by an algebraic topology argument.\footnote{One can consult the paper of Akbulut and Larson \cite{AL18} for the handle diagram of a rational homology $4$-ball including a $3$-handle. This $4$-manifold has the boundary $\Sigma(2,3,7)$.} Then, the above characterization indicates that such a rational homology $4$-ball cannot be a Stein domain. Mazur manifolds are potential candidates of Stein domains, but this is not the case for all Mazur manifolds, see the impressive work of Mark and Tosun \cite{MT18}.

In addition to the non-injectivity of $\psi$, we know that it is not surjective. In particular, Kim and Livingston proved that $\mathrm{Coker}(\psi)$ has a $\Z^\infty \oplus \Z_2^\infty$ subgroup \cite{KS14}. This was reproved by Aceto and Larson \cite{AL17} as a consequence of a more general fact. They proved that $\psi \left( \TZ \right)$ and $\mathcal{L}$ intersect trivially where $\mathcal{L}$ denotes the subgroup of $\TQ$ generated by lens spaces. In particular, the structure of $\mathcal{L}$ has been studied in \cite{AL17, ACP20}. Finally, we can ask:

\begin{introprob}
Does $\mathrm{Coker}(\psi)$ contain a $\Z^\infty \oplus \Z_2^\infty$ summand? Does it have other types of subgroups or summands?
\end{introprob}

In the light of the results therein and in Section~\ref{trivial}, we can also address the following explicit problem:

\begin{introprob}
Do the Brieskorn spheres $\Sigma(2,3,6n+1)$ bound rational homology $4$-balls (resp. homology $4$-balls) for odd $n \geq 5$ (resp. even $n \geq 6$)?
\end{introprob}

The notion of rational homology cobordism can be generalized among all closed connected oriented $3$-manifolds. Such a homology cobordism is said to be \emph{ribbon} if the cobordism $4$-manifold is built by attaching only $1$- and $2$-handles. This give rises to a preorder on the set of homeomorphism classes of closed connected oriented $3$-manifolds. Daemi, Lidman, Vela-Vick, and Wong conjectured that this preorder is in fact a partial order. Inpendently, Friedl, Misev, and Zentner \cite{FMZ22} and Huber \cite{H22} proved this conjecture affirmatively, relying on the result of Agol \cite{A22}.

\section{Appendix: Examples of Homology \texorpdfstring{$3$}{3}-Spheres}

\label{appendix}

In the wide world of closed connected oriented $3$-manifolds, there is a simple characterization of homology $3$-spheres $Y$ thanks to Poincar\'e duality and universal coefficient theorem: $H_1(Y;\Z) = 0$. Since the abelianization of $\pi_1 (Y)$ gives $H_1(Y;\Z)$ due to Hurewicz theorem, they are even easily recognized. In this appendix, we discuss several constructions of homology $3$-spheres, our main references are Neumann and Raymond \cite{NR78}, Eisenbud and Neumann \cite{EN85}, Gompf and Stipsicz \cite{GS99}, Saveliev \cite{Sav02}, and Akbulut \cite{A16}. 

The first example of homology $3$ spheres was given by Poincar\'e \cite{P04} as a counterexample to the first version of the Poincar\'e conjecture. This $3$ manifold is known as \emph{Poincare homology sphere} and the exposition of Kirby and Scharlemann can be seen for the eight equivalent descriptions of the Poincar\'e homology sphere \cite{KS79}.

The next source for homology $3$-spheres was found by Dehn \cite{D38} by providing a passage from $1$-manifolds -knots and links- to $3$-manifolds via the topological operation called \emph{surgery}. Consider the tubular neighborhood of $K$ in $S^3$, which is a solid torus $\nu(K) \approx S^1 \times D^2$. On the boundary torus $\partial \nu(K)$, there is a preferred longitude $\lambda$, i.e., a simple closed curve with $\mathrm{lk}(\lambda,K) = 0$, and there is a canonical meridian $\mu$ with $\mathrm{lk}(\mu,K) = 1$. 

A \emph{Dehn $(p/q)$-surgery} along $K$ in $S^3$ is constructed by following two steps. We first drill out the interior of $\nu(K)$ from $S^3$ and consider the knot exterior $S^3 \setminus \intB$. Next, we glue another solid torus $D^2 \times S^1$ to the knot exterior by a homeomorphism $\varphi$. The resulting closed $3$-manifold $S^3_{p/q}(K)$ is given by \begin{center}
$S^3_{p/q}(K) = \left (S^3 \setminus \intB \right) \cup_\varphi \left( D^2 \times S^1 \right), \ \ \ \varphi(\partial D^2 \times \{ *\}) = p\mu + q\lambda.$ \end{center} Since $H_1(S^3_{p/q}(K); \Z) = \Z_p$, the manifolds of the form $S^3_{1/n}(K)$ are automatically homology $3$-spheres. In particular, Dehn showed that the Poincar\'e homology sphere can be obtained by $(-1)$-surgery along the left-handed trefoil knot $T(2,3)$ in $S^3$.

A \emph{framed knot} in $S^3$ is a knot equipped with a smooth nowhere vanishing vector field normal to the knot. Thus a \emph{framing} of a knot is naturally characterized by its Seifert surface \cite{S35} and \cite{FP30} so that the specified longitude is given by $0$-framing.\footnote{The existence of Seifert surfaces of an oriented knot $K$ in an oriented $3$-manifold $M$ would be possible if and only if $K$ is null-homologous, i.e., $[K] = 0 \in H_1(M;\Z)$, one can consult \cite{R76}.} The set of framings of a knot is identified with a fixed set of rationals using a Seifert surface, so each knot has a preferred well-defined framing. This process can be naturally generalized to framed \emph{links} in $S^3$, which are disjoint collections of knots in $S^3$.

By the eminent results of Lickorish \cite{L62} and Wallace \cite{W60}, and Kirby \cite{K78}: the map $D$ provided by integral $n$-surgery \begin{center}
$D: \{ \text{framed links in} \ S^3 \} \to \{ \text{closed} \ 3\text{-manifolds}\}, \ \ \ L \mapsto D(L) = S^3_n (L) $
\end{center} is many-to-one. In particular, Kirby completely described when two elements can represent the same element in the kernel using Cerf theory \cite{C70}, i.e., $S^3_n (L_1)$ is homeomorphic to  $S^3_n (L_2)$ if and only if the framed links are related by sequences of two \emph{Kirby moves}: blow-up and handle-slide. His notable contribution was generalized, ramified, and reproved by Fenn and Rourke \cite{FR79}, C\'{e}sar de S\'{a} \cite{C79}, Kaplan \cite{K79}, Rolfsen \cite{R84}, Lu \cite{L92}, Matveev and Polyak \cite{MP94}, and Martelli \cite{Ma12}.

The next construction of homology $3$-spheres was provided by Seifert \cite{S33}. Let $e$ be an integer and let $(a_1,b_1), \ldots, (a_n,b_n)$ be pairs of relatively prime integers. The \emph{Seifert fibered space} with base orbifold $S^2$ is a closed $3$-manifold \begin{center}
$M(S^2;e, (a_1,b_1), \ldots, (a_n,b_n))$
\end{center} constructed by starting with $S^1$-bundle over an $n$-punctured $S^2$ of Euler number $e$ and filling the $k$th boundary component with $\left ( a_k/b_k \right )$-framed solid torus for $k=1, \ldots, n$. The core circle of the $\left ( a_k/b_k \right )$ Dehn filling is called a \emph{singular fiber}, all other fibers are said to be \emph{regular fibers}. The resulting manifold is a homology $3$-sphere if and only if  \begin{align}
a_1 \ldots a_n \left ( -e + \sum^n_{k=1} \frac{b_k}{a_k} \right )  = \mp 1.
\label{weights} 
\end{align} This equation results from the fundamental group \cite[Pg. 398]{ST80}, and hence the first homology group calculations of Seifert fibered spaces, see \cite[Pg. 410]{ST80}.\footnote{Seifert called homology $3$-spheres \emph{Poincar\'e spaces}, see \cite[Pg. 402]{ST80}. Note that the book \cite{ST80} includes an English translation of \cite{S33} and our citations all lie in that part.} In particular, Poincar\'e homology sphere corresponds to the Seifert fibered space $M(S^2; -2, (2,-1),(3,-2),(5,-4))$.

Due to Brieskorn \cite{B66a, B66b}, homology $3$-spheres also originate from algebraic geometry as seen in the variety of certain complex analytical polynomials. Let $p,q$ and $r$ be relatively coprime positive integers. Let $f: \mathbb{C}^3 \to \mathbb{C}$ be a complex analytical polynomial defined by $f(x,y,z) = x^p + y^q + z^r$. Then the zero set of $f$ is the complex surface $V(f) = \{ (x,y,z) \in \mathbb{C}^3 \ \vert \ f(x,y,z)=0 \}$ singular at the origin. If we transversally intersect this variety with the five-sphere $S^5_\epsilon$ of arbitrarily small radius $\epsilon$, then the resulting closed $3$-manifold is the \emph{Brieskorn sphere} given by \begin{center} $\Sigma(p,q,r) = V(f) \pitchfork S^5_\epsilon \subset \mathbb{C}^3 .$ \end{center} The Poincar\'e homology sphere matches with the Brieskorn sphere $\Sigma(2,3,5)$. For explicit descriptions of fundamental groups of Brieskorn spheres, see Milnor's paper \cite{M75}. In particular, there is an orientation-preserving homemorphism between $M(S^2;a_1,a_2,a_3)$ and $\Sigma(a_1,a_2,a_3)$ \cite[Theorem~4.1]{NR78}. In general, it is possible to realize Seifert fibered homology $3$-spheres as the links of the complex surface singularities of Brieskorn complete intersections \begin{center}$
V_B (a_1, \ldots, a_n) = \{ b_{i1} z_1^{a_1} + \ldots + b_{in} z_n^{a_n} =0, \ i=1, \ldots, n-2  \} \subset \mathbb{C}^n$
\end{center} where $B = (b_{ij})$ is an $(n-2) \times n$-matrix of complex numbers such that each of the maximal minors of $B$ is non-zero, see \cite[Theorem~2.1]{NR78}.

Let $\mathcal{J}$ be an index set. A \emph{plumbing graph} $G$ is a connected and weighted tree with vertices $v_j$ and weights $e_j$ for $j \in \mathcal{J}$. We can construct a $4$-manifold $X(G)$ with a boundary $Y(G)$ by using the plumbing graph. First, for each $v_j$, we assign a $D^2$-bundle over $S^2$ whose Euler number is $e_j$. Next, we plumb two of these $D^2$-bundles if there is an edge connecting the vertices, see \cite[Theorem~5.1]{NR78}.  

The fundamental classes of the zero-sections of $D^2$-bundles generate the second homology group $H_2(X(G); \mathbb{Z})$. Thus, for each vertex of $G$, we have a generator of $H_2(X(G); \mathbb{Z})$. Hence, the intersection form on $H_2(X(G); \mathbb{Z})$ is naturally characterized by the corresponding intersection matrix $I=(a_{ij})$ whose data is given in the following way: \[ a_{ij} = \begin{cases} 
      e_i, & \text{if} \ v_i=v_j, \\
      1, & \text{if} \ v_i \ \text{and} \ v_j \ \text{is connected by one edge}, \\
      0, & \text{otherwise}. 
\end{cases}\]
   
A plumbing graph $G$ is called \emph{unimodular} if $\mathrm{det}(I)= \pm 1$. The unimodularity of the plumbing graph implies that $Y(G)$ is a homology $3$-sphere, so it is called a \emph{plumbed homology $3$-sphere}. We may characterize the \emph{negative definiteness} of $G$, it requires that $I$ is negative-definite, i.e., $\mathrm{signature}(I)=-|G|$, where $|G|$ denotes the number of vertices of $G$. 

A Seifert fibered homology sphere $M(S^2;e, (a_1,b_1), \ldots, (a_n,b_n))$ can be realized as the boundary of a star-shaped plumbing graph. This graph is unique when it is negative-definite \cite[Section~1.1]{Sav02}. The integer weights $t_{ij}$ in the graph are found by solving the equation \eqref{weights} and expanding the continued fractions $[t_{i1},\ldots,t_{im_{i}}]$ as follows: for each $i \in \{1, \ldots,n \}$, we have

\vspace{0.3 em} 

\begin{minipage}{0.6\linewidth}

\centering

\[ a_i/b_i = [t_{i1},t_{i2}, \ldots,t_{im_{i}}]=t_{i1}-\cfrac{1}{t_{i2}-\cfrac{1}{\cdots-\cfrac{1}{t_{im_{i}}}}}  \]
  
\end{minipage}  
\begin{minipage}{0.4\linewidth}  
\centering
\includegraphics[width=0.5\textwidth]{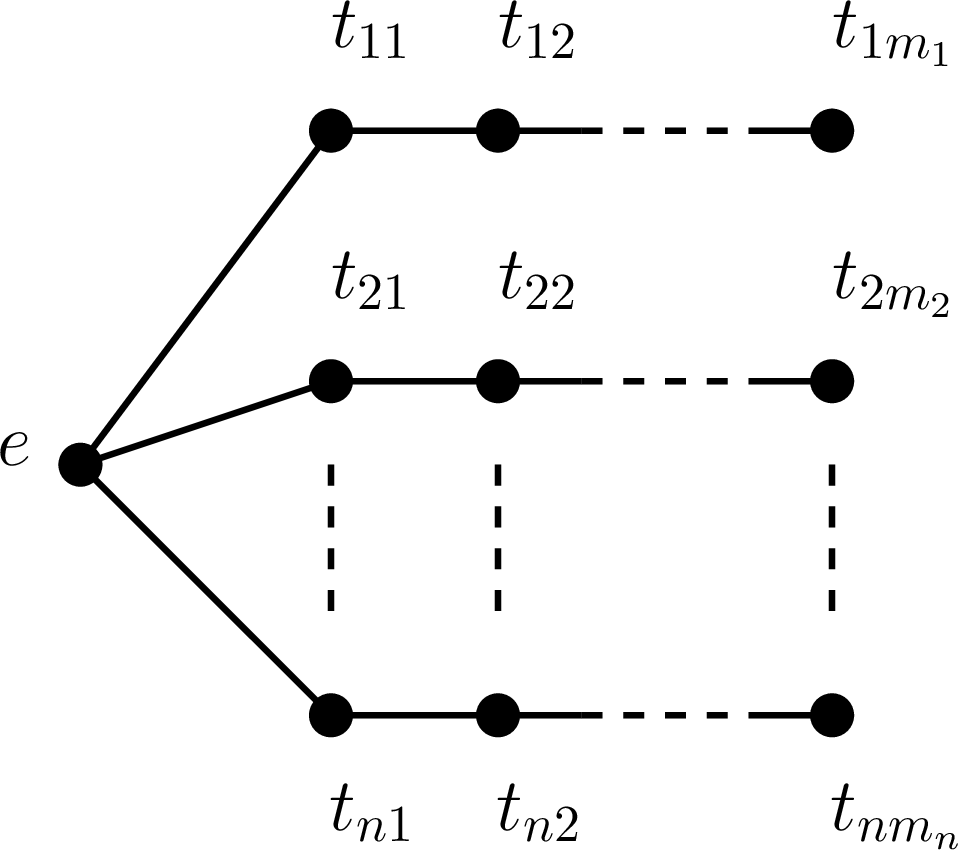}   
\end{minipage}  

\vspace{0.7 em} 
 
In this survey, we focus on the following three families of Brieskorn spheres. Assume that $p$ and $q$ are pairwise coprime, positive, and ordered integers such that $2 \leq p < q$:

\begin{enumerate}
\item \label{1} $\{ \Sigma(p,q,pqn-1)\}_{n=1}^\infty$,
\item \label{2} $\{ \Sigma(p,q,pqn+1)\}_{n=1}^\infty$,
\item \label{3} $\{ \Sigma(p_n, q_n, r_n)\}_{n=1}^\infty$ where $p_n q_n + p_n r_n - q_n r_n = 1$,
\begin{itemize}
\item[(a)] $\{ \Sigma(2n, 4n-1,4n+1)\}_{n=1}^\infty$,
\item[(b)] $\{ \Sigma(2n+1, 4n+1,4n+3)\}_{n=1}^\infty$,
\item[(c)] $\{ \Sigma(2n+1, 3n+2,6n+1)\}_{n=1}^\infty$,
\item[(d)] $\{ \Sigma(2n+1, 3n+1,6n+5)\}_{n=1}^\infty$.
\end{itemize}
\end{enumerate}

Due to the classical result of Moser \cite{M71}, the first two families can be obtained by $(-1/n)$ surgeries along the left-handed torus knots $T(p,q)$ and their mirrors right-handed torus knots $\overline{T(p,q)}$ in $S^3$: \begin{center} $\Sigma(p,q,pqn-1) = S^3_{-1/n}(T(p,q)), \ \text{and} \ \Sigma(p,q,pqn+1) = S^3_{-1/n}(\overline{T(p,q)}).$ \end{center} The third family is called \emph{almost simple linear graphs} and extensively studied in \cite{FS85}, \cite{E95}, and \cite{KS20}. The families \eqref{1} and \eqref{3} are vast generalizations of the Poincar\'e homology sphere $\Sigma(2,3,5)$ while the the family \eqref{2} is of $\Sigma(2,3,7)$. 

Note that there is a family of Brieskorn spheres realized as the boundaries of almost simple graphs which cannot be obtained by surgeries along any knots in $S^3$. This surgery obstruction was due to Hom, Karakurt, and Lidman \cite{HKL16}. In particular, they showed that $\Sigma(2n, 4n-1,4n+1)$ cannot be realized as knot surgeries for $n \geq 4$. 

Another classical way to produce homology $3$-spheres is the method of \emph{cyclic branched coverings} of $S^3$ branched over knots $K$, which dates back to work of Alexander \cite{A20} and Seifert \cite{S33}. Let $X_n (K)$ be the $n$-fold regular covering of the knot exterior $ X(K) = S^3 \setminus \intB$. Then the $n$-fold cyclic branched covering of $S^3$ over $K$ is a closed $3$-manifold \begin{center} $\Sigma_n (K) = X_n (K) \cup_\varphi \left( D^2 \times S^1 \right), \ \ \ \varphi( \tilde{\mu} )=  \mu$ \end{center} where $\mu \subset \partial X (K)$ is the meridian of $K$ and $\tilde{\mu}$ is the lift of $\mu$ to $\partial X_n (K)$. Note that $\Sigma_n (K)$ is a homology $3$-sphere when \[ \prod_{k=1}^n \Delta_K \left ( e^{\frac{2\pi i k}{n}} \right ) = 1 \] where $\Delta_K (t)$ is the Alexander polynomial of $K$ normalized so that there are no negative powers of $t$ and the constant term is positive. The Brieskorn sphere $\Sigma(p,q,r)$ is an $r$-fold cyclic branched coverings of $S^3$ branched over the torus knots $T(p,q)$, see \cite[Lemma~1.1]{M75} and compare with \cite[Pg. 405]{ST80}. In general, Seifert fibered sphere $\Sigma(a_1, \ldots, a_n)$ is a $2$-fold cyclic branched coverings of $S^3$ branched over Montesinos knots $K(a_1, \ldots, a_n)$, see \cite{Mon73, Mon75}.

Given two homology $3$-spheres together with knots inside them, we can produce a new closed $3$-manifold by following the agenda of Gordon \cite{G75}. 

Let $K_1$ and $K_2$ be knots in homology $3$-spheres $Y_1$ and $Y_2$ with the knot exteriors $Y_1 \setminus \intC$ and $Y_2 \setminus \intD$, and the longitude-meridian pairs $(\lambda_1, \mu_1)$ and $(\lambda_2, \mu_2)$ respectively. Consider the following integral $2 \times 2$ matrix $A = \begin{pmatrix} a & b \\ c & d \end{pmatrix}$ with $\mathrm{det}(A) = -1.$ Gordon constructed closed $3$-manifolds obtained by gluing knot exteriors of homology $3$-spheres along their boundary tori by matching longitude-meridian pairs with respect to the matrix $A$: \begin{center}
$Y(K_1,K_2,A) = \left ( Y_1 \setminus \intC \right) \cup_A \left ( Y_2 \setminus \intD \right) .$ \end{center} Clearly, the resulting manifold is a homology $3$-sphere whenever $A = \begin{pmatrix}  a & ab+1 \\ 1 & b \end{pmatrix}$. Gordon studied the problem which $Y(K_1,K_2,A)$ bound contractible $4$-manifolds and provided several characterizations in terms of sliceness of knots.
\vspace*{0.4 em}

The case $A = \begin{pmatrix} 0 & 1 \\ 1 & 0 \end{pmatrix}$ corresponds to switching longitude-meridian pairs of knots inside homology $3$-spheres. This construction is of special interest and is known as the \emph{splice} operation first introduced by Siebenmann \cite{S80}. Given the pairs $(Y_1,K_1)$ and $(Y_2,K_2)$, we will denote the splice of these manifolds along the given knots by $Y_1 \tensor[_{K_1}]{\bowtie}{_{K_2}} Y_2.$

The concept of the splice became popular after the novel book of Eisenbud and Neumann \cite{EN85} because the splice can be realized as a generalization of several other topological operations including cabling, connected sum, and disjoint union. The splice also has a very crucial role in singularity theory due to Neumann and Wahl \cite{NW90}. For details, one can consult the recent survey of Cueto, Popescu-Pampu, and Stepanov \cite{CPPS22}.

We finally consider the graph $3$-manifolds introduced by Waldhausen \cite{W67}. A \emph{graph $3$-manifold} is a closed $3$-manifold such that it can be cut along a set of disjoint embedded tori $T_i$ and has a decomposition with each piece is $\Sigma_i \times S^1$, where $\Sigma_i$ is a surface with boundary. In the light of JSJ (torus) decomposition theorem (Jaco and Shalen \cite{JS79} and Johannson \cite{J79}), a \emph{graph homology $3$-sphere} is a prime homology $3$-sphere whose JSJ decomposition contains only Seifert fibered pieces. See Neumann's paper \cite{N07} and its appendix, and Saveliev's book \cite{Sav02} for further discussions.

\printendnotes

\section*{Afterword}

The recorded history of the $n$-dimensional homology cobordism group $\HCZ$ first appeared in the Ph.D. thesis of Gonz\'{a}lez-Acu\~{n}a \cite{GA70} under the supervision of Ralph H. Fox at Princeton University in 1970. He introduced this notion to study homology $n$-spheres by building on the work of Kervaire and Milnor \cite{KM63} about the $n$-dimensional homotopy cobordism group $\HC$ of homotopy $n$-spheres. Gonz\'{a}lez-Acu\~{n}a proved that these groups $\HC$ and $\HCZ$ are isomorphic unless $n=3$. Therefore, they are both finite except in the case of $n=3$. This result was not published as an article but was referred in \cite[Section~2]{G70}. Note that the only unknown value of the order of $\HC$ in \cite{KM63} was the case of $n=3$. This has not been clarified until the work of Perelman \cite{P1, P3, P2}.

The isomorphism argument of Gonz\'{a}lez-Acu\~{n}a broke down when $n=3$, if the order of $\Theta^3$ was known at that time, see \cite[Pg.~17, Remark and Section~I.5]{GA70}. Especially, the homology cobordism group $\TZ$ was introduced to him by Denis Sullivan as noted in \cite[Pg.~VII]{GA70}.  Also, the first known proof of the homology cobordism invariance of the Rokhlin invariant $\mu$ was given \cite[Pg.~33-34]{GA70}. Further, the relation between the Arf invariant of knots and the Rokhlin invariant in terms of  knot surgery was found \cite[Theorem~III.2]{GA70}. Unfortunately, his results were only mentioned in Gordon's article \cite{G75} and they have remained mysteries. 

The main references for our survey are the great book of Saveliev \cite{Sav02} and the eminent ICM 2018 article of Manolescu \cite{M18}. To extend their coherent frameworks, we list recent results not included in these resources. Further, we catalog all natural sources of homology $3$-spheres in the appendix. In order to avoid distracting readers, we share our footnotes as endnotes. 

\section*{Acknowledgements}

The author would like to thank Selman Akbulut, Ronald Fintushel, Yoshihiro Fukumoto, Kristen Hendricks, Jennifer Hom, \c{C}a\u{g}r{\i} Karakurt, Nikolai Saveliev, Steven Sivek, Ronald Stern, Andr\'as Stipsicz, and Zhouli Xu for their valuable comments and feedback on earlier drafts of this document. Special thanks go to Tye Lidman, Ciprian Manolescu, and Masaki Taniguchi for sharing their expertise on the subject and providing insightful suggestions. Otherwise, the survey would have been incomplete from several perspectives. Also, we are grateful to Francisco Javier Gonz\'{a}lez-Acu\~{n}a for sharing a scanned copy of his article \cite{G70}. 

This survey has been conducted respectively at Max-Planck-Institut f\"ur Mathematik in Bonn, Bo\u{g}azi\c{c}i University in Istanbul, and Stanford University in California. We are indebted to all these institutions for their generous hospitality and support. The author is currently supported by the Turkish Fulbright Commission \say{Ph.D. dissertation research grant}.

Finally, the author would like to thank the anonymous referee for the invaluable feedback that improved both the content and the exposition of the survey.

\bibliography{homologycobordismgroup}
\bibliographystyle{amsalpha}

\end{document}